\documentclass[leqno,11pt]{article}
\usepackage{latexsym}
\usepackage{amssymb}
\usepackage{amsfonts}
\usepackage{amsmath}
\usepackage{color}

\allowdisplaybreaks
\makeatletter
\long\def\unmarkedfootnote#1{{\long\def\@makefntext##1{##1}\footnotetext{#1}}}
\makeatother

\newtheorem{definition}{Definition}[section]

\newtheorem{lemma}[definition]{Lemma}

\newtheorem{theorem}[definition]{Theorem}

\newtheorem{proposition}[definition]{Proposition}

\newtheorem{corollary}[definition]{Corollary}

\newtheorem{remark}[definition]{Remark}

\newtheorem{example}[definition]{Example}

\newtheorem{examples}[definition]{Examples}

\def\o{\Omega}

\def\w0{{W_0^{m,A}(\Omega)}}

\def\R{\mathbb R}
\def\N{\mathbb N}
\def\ep{\varepsilon}
\newcommand{\medint}{-\kern  -,395cm\int}
\newcommand{\medintinrigo}{-\kern  -,315cm\int}
\newcommand{\medelle}{-\kern  -,235cm L}
\newcommand{\medellenrigo}{-\kern  -,180cm L}
\newcommand{\qed}{\thinspace\null\nobreak\hfill
\hbox{\vbox{\kern-.2pt\hrule height.2pt
depth.2pt\kern-.2pt\kern-.2pt \hbox to1.8mm {\kern-.2pt\vrule
width.4pt \kern-.2pt\raise1.8mm\vbox to.2pt{} \lower0pt\vtop
to.2pt{}\hfil\kern-.2pt \vrule
width.4pt\kern-.2pt}\kern-.2pt\kern-.2pt \hrule height.2pt
depth.2pt \kern-.2pt}}\par\medbreak}

\newcommand{\rn}{\mathbb R^n}

\newcommand{\np}{{n^\prime}}

\newcommand{\LLN}{{\cal L}^n}

\newcommand{\HN}{{\cal H}^{n-1}}
\newcommand{\HD}{{\cal H}^{d}}

\newcommand{\cn}{n \omega _n^{1/n}}
\newcommand{\on}{\omega _n}

\newcommand{\wpp}{W^{1,p}(\o )}

\newcommand{\bro}{B_\rho (x_0)}
\newcommand{\wnn}{W^{1,n}(\o )}
\newcommand{\wno}{W^{1,n}_0(\o )}
\newcommand{\wln}{W^1L^{n,q}(\o )}
\newcommand{\wlo}{W^1_0L^{n,q}(\o )}
\newcommand{\mo}{|\o |}

\title{Continuity properties of weakly monotone Orlicz-Sobolev functions} 
\numberwithin{equation}{section}
\author{
Menita Carozza\\
{\it  Dipartimento di Ingegneria, Universit\`a del Sannio}
\\ {\it Corso Garibaldi 107, 82100 Benevento, Italy}
\\ {\it e-mail: carozza@unisannio.it }
\bigskip
\\
  Andrea Cianchi\\
 {\it Dipartimento di Matematica e Informatica \lq\lq U. Dini", Universit\`a di Firenze}\\ {\it Viale Morgagni 67/A, 50134 Firenze, Italy}
 \\{\it  e-mail: cianchi@unifi.it}}

\begin{document}
\maketitle

\begin{abstract}
The notion of weakly monotone functions extends the classical definition of monotone function, that can be traced back to H.Lebesgue. It was introduced, in the setting of Sobolev spaces, by J.Manfredi, and  thoroughly investigated in the more general framework of Orlicz-Sobolev spaces by diverse authors, including T.Iwaniec, J.Kauhanen, P.Koskela, J.Maly, J.Onninen, X.Zhong.
The present paper complements and augments  the available theory of pointwise regularity properties of weakly monotone functions in Orlicz-Sobolev spaces. In particular, a variant is proposed in a customary condition ensuring the continuity of functions from these spaces which avoids a technical additional assumption, and applies to certain situations when the latter is not fulfilled. The continuity outside sets of zero Orlicz capacity, and outside sets of (generalized) zero Hausdorff measure, will are also established  when everywhere continuity fails.
\end{abstract}

\unmarkedfootnote {
}

\maketitle

\section{Introduction}\label{intro}

A weakly monotone function in an open set $\o \subset \rn$, $n \geq 2$, is, loosely speaking, a Sobolev function that satisfies the minimum and maximum principles in a weak sense. Precisely, a function $u \in  W_{\rm loc}^{1,1}(\Omega)$ is called weakly monotone if, 
 for every open set $\o ' \subset \subset\Omega$,   and   every  $m,\, M \in \mathbb R$, such that $m\le M$ and
$$(u-M)_+\,-\,(m-u)_+\in W_0^{1,1}(\o')\,,$$
one has  that
$$ m\le u\le M\qquad  \hbox{a.e.  in $\o'$.}$$
Here, the subscript $+$ stands for positive part. 
\par The notion of weak monotonicity was introduced by Manfredi in \cite{M}, where he provided  a new direct approach to the regularity theory of maps with finite distortion, and of maps in classes defined in terms of integrability properties of the adjugates of their gradients, which play a role in nonlinear elasticity (see \cite{B}). Earlier proofs of continuity properties of these maps, contained in  \cite{GV} and   \cite{S},  made use of the notion of topological degree. 
\par A key idea in \cite{M} is to exploit the fact that the components of these maps are weakly monotone functions, and that any such function is continuous, or at least continuous outside a set of a certain capacity zero, provided that a sufficiently large power of the modulus of its gradient is integrable.  Specifically, assume that $u$ is a weakly  monotone function from $W^{1,p}_{\rm loc}(\o)$ for some $p \geq 1$. If $p>n$, then $u$ is continuous (irrespective of whether it is weakly monotone or not), by the Sobolev embedding theorem. Hence, it is monotone in the classical sense introduced by Lebesgue in his study of the Dirichlet problem in the plane \cite{L}. The advance of \cite[Theorem 1]{M} amounts to showing that, even
\begin{equation} \label{p=n}   \hbox{if $p=n$, then $u$ is continuous,}
\end{equation}
and that
\begin{equation} \label{p>n-1} \hbox{if $p>n-1$, then  $u$  is continuous outside  a set of $C_{p,1}$-capacity zero.}
\end{equation}
\par Weakly monotone functions  come into play in the regularity theory of elliptic partial differential equations as well. For instance, as pointed out in \cite{HKM} and \cite{KMV}, weak solutions to  $p$-Laplacian type  elliptic equations, with possibly degenerating ellipticity, turn out to be weakly monotone.
\par
The result of \cite{M} has paved the way to investigations on pointwise  properties of weakly monotone functions in more general classes of Sobolev type spaces. In particular, weakly monotone functions from Orlicz-Sobolev spaces $W^{1,A}_{\rm loc}(\o)$ are focused in the monograph \cite{IM}, and in the papers \cite{IKO, KKMOZ}. A motivation for these studies was  the analysis of maps of bounded distortion whose gradient locally belongs to the  local Orlicz space $L^A_{\rm loc}(\o)$ for some  Young function $A$ that is not necessarily of power type. These contributions pointed out that continuity of a weakly monotone function is guaranteed even if it belongs to an Orlicz-Sobolev space slightly larger than $W^{1,n}_{\rm loc}(\o)$, namely if $A(t)$ grows slightly more slowly than $t^n$ near infinity. Precisely,  \cite[Proposition 2.7]{KKMOZ} states that, if
\begin{equation}
\label{Iwcond}
\int ^\infty \frac{A(t)}{t^{n+1}}\, dt =
\infty \,,
\end{equation}
and there exists $\varepsilon >0$ such that the function
\begin{equation}
\label{incr}
 t \mapsto \frac{A(t)}{t^{n-1+\varepsilon}} \qquad \hbox{is increasing,}
\end{equation}
then any weakly monotone function from $W^{1,A}_{\rm loc}(\o)$ is continuous.  The same conclusion, with \eqref{incr} replaced by a  slightly stronger condition of a similar nature,  is proved in \cite[Theorem 7.5.1]{IM}. Furthermore, information on its (local) modulus of continuity is provided. Assumptions of a different kind for the continuity of weakly monotone functions in Orlicz-Sobolev spaces can also be found in \cite{FM}.
\par
A condition of the form \eqref{Iwcond} amounts to imposing an appropriate degree of integrability of  the gradient of the weakly monotone functions in question, and is an indispensable   requirement. On the other hand, assumption \eqref{incr} has an essentially technical nature. In fact, a close inspection of the proof of \cite{KKMOZ} reveals that \eqref{incr} is basically needed to deduce certain properties of Orlicz-Sobolev functions from their analogoues in the  theory of standard Sobolev spaces.
\par
In the present paper, we suggest some variants in the approach of \cite{M}, \cite{IM} and \cite{KKMOZ}, that call into play peculiar Orlicz space techniques and results. This enables us to drop condition \eqref{incr}, and to establish the everywhere continuity of weakly monotone functions from the space  $W^{1,A}_{\rm loc}(\o)$ under a single assumption, in the spirit of \eqref{Iwcond}, but with $A$ replaced by a closely related Young function depending also on $n$, that will be denoted by $A_{n-1}$. Namely, our condition reads
\begin{equation}
\label{ours}
\int ^\infty \frac{A_{n-1}(t)}{t^{n+1}}\, dt =
\infty \,,
\end{equation}
and it also implies the local uniform continuity of any weakly monotone function in $W^{1,A}_{\rm loc}(\o)$, with an explicit modulus of continuity depending only on $A$ and $n$.
The function $A_{n-1}$ comes into play in a sharp Poincar\'e type inequality for the oscillation of Sobolev functions on the $(n-1)$-dimensional unit sphere in $\rn$. A definition of $A_{n-1}$ can be found in Section \ref{main}, where the main results of this paper are stated. Here, let us just mention  that, if $n=2$, assumption \eqref{ours} 
coincides with \eqref{Iwcond}, since $A_1=A$. This shows that condition \eqref{incr} is actually irrelevant in the results of \cite{IM, KKMOZ} in this case. When $n \geq 3$, the function $A_{n-1}$  is  equivalent to $A$, and hence conditions \eqref{Iwcond} and \eqref{ours} again agree, in any customary, non-borderline situation. Here, loosely speaking, borderline means that $A(t)$  is not larger than  $t^{n-1}$ near infinity, in which case the function $A_{n-1}$ can grow slightly faster that  $A$ near infinity.
\par
This is the content of Theorem \ref{febbraio0}, which enhances the results of \cite{IM, KKMOZ}, since  conditions \eqref{Iwcond}--\eqref{incr} imply \eqref{ours}, 
 whereas Young functions can be exhibited that fulfill \eqref{ours}, but not \eqref{incr} -- see Proposition \ref{improve}. This shows that Theorem \ref{febbraio0} is applicable in circumstances where the  available results in the literature may fail. Moreover,  the results of \cite{IM, KKMOZ} can be recovered  as a consequence of Theorem \ref{febbraio0}, and, in fact, condition \eqref{Iwcond} can be shown to be sufficient for the continuity of weakly monotone functions from $W^{1,A}_{\rm loc}(\Omega)$ with an additional condition slightly less demanding  than \eqref{incr} -- see Corollary \ref{recover}.
\par Under a weaker assumption than \eqref{ours} -- a counterpart of the assumption $p>n-1$ appearing in \eqref{p>n-1} for classical Sobolev spaces -- 
in Theorem \ref{orliczbound} we  prove that every weakly monotone function from $W^{1,A}_{\rm loc}(\o)$ is locally bounded and differentiable a.e. in $\o$. The assumption in question is only needed when $n\geq 3$, and takes the form
\begin{equation}\label{equiv2bis} \int ^\infty \bigg(\frac{t}{A(t)}\bigg)^{\frac 1{n-2}} \, dt <
\infty\,.
 \end{equation}
If $n=2$, the same conclusion holds whatever $A$ is.
\par
Under the same  hypothesis on $n$ and $A$, every weakly monotone function from $W^{1,A}_{\rm loc}(\o)$  is shown to be continuous outside an exceptional set of vanishing Orlicz capacity, that depends on the function $A$. This is the subject of Theorem \ref{orliczcapacity}, that not only extends, but also somewhat augments property \eqref{p>n-1} even in the case when  $W^{1,A}_{\rm loc}(\Omega)=  W^{1,p}_{\rm loc}(\Omega)$. This is observed in Remark \ref{improveManfredi}. Having Theorem \ref{orliczcapacity} at disposal, an estimate for the size of the singular set in terms of a Hausdorff measure,  defined in terms of $A$, is established   in  Theorem \ref{hausdorff}.

\section{Orlicz and Orlicz-Sobolev spaces}\label{youngorlicz}

\color{black}
%
%

The notion of Orlicz space relies upon that of Young function. 
 A  function  $A: [0,
\infty ) \to [0, \infty ]$ is called a Young function if it is
convex, non constant in $(0, \infty)$, and  vanishes at $0$. Any function fulfilling these properties
has the form
\begin{equation}\label{A}
A(t) = \int _0^t a(r )\, dr \qquad \quad \hbox{for $t \geq 0$},
\end{equation}
for some non-decreasing, left-continuous function $a: [0, \infty )
\to [0, \infty ]$ which is neither identically $0$, nor infinity. Observe that the function
\begin{equation}\label{monotone}
\hbox{$\displaystyle t \mapsto \frac{A(t)}{t}$ \quad is
non-decreasing,}
\end{equation}
and 
\begin{equation}\label{gen9}
A(t) \leq a(t)\,t \leq A(2t) \quad \hbox{for $t\geq 0$.}
\end{equation}
Furthermore,  if $k \geq 1$, then
\begin{equation}\label{kt}
kA(t) \leq A(kt) \quad \hbox{for $t \geq0$,}
\end{equation}
and hence
\begin{equation}\label{kt-1}
kA^{-1}(t) \geq A^{-1}(kt) \quad \hbox{for $t \geq0$.}
\end{equation}
Here, $A^{-1}$ denotes the (generalized) right-continuous inverse of $A$.  The Young conjugate $\widetilde{A}$ of $A$  is
defined by
$$\widetilde{A}(t) = \sup \{s t-A(s):\,s \geq 0\} \qquad {\rm for}\qquad  t\geq 0\,.$$
The alternative notation $A^{\widetilde {}}$ 
will also be adopted instead of for $\widetilde {A}$ whenever convenient. Note
the representation  formula
\begin{equation}\label{youngconj}
\widetilde A(t) = \int _0^t a^{-1}(r ) \,dr \qquad \quad
\hbox{for $t \geq 0$,}
\end{equation}
where $a^{-1}$ denotes the (generalized) left-continuous inverse of
the function $a$. Let us notice that $\widetilde {\widetilde A}=A$.
\\ If, for instance, $A(t) = \tfrac {t^p}p$ for some $p \in (1, \infty)$, then $\widetilde A (t) =\tfrac {t^{p'}}{p'}$, where $p'= \tfrac{p}{p-1}$, the H\"older conjugate of $p$. 
\\ A property to be used in what follows is that, if $A$ is a Young function and $q \in (1, \infty)$, then 
\begin{equation}\label{equivconj}
\hbox{the function $t \mapsto \displaystyle \frac {A(t)}{t^q}$\,\, is increasing if and only if the function \, $t\mapsto\displaystyle \frac {\widetilde A(t)}{t^{q'}}$\,\, is decreasing\,.}
\end{equation}
An application of equation
\eqref{gen9} with $A$ replaced by $\widetilde A$ yields
\begin{equation}\label{gen10}
\widetilde A(t) \leq a^{-1}(t)\,t \leq \widetilde A(2t) \quad
\hbox{for $t\geq 0$.}
\end{equation}
Moreover, one has that
\begin{equation}\label{AAtilde}
t \leq A^{-1}(t) \widetilde A^{-1}(t) \leq 2t \qquad \hbox{for $ t
\geq 0.$}
\end{equation}
A Young function $A$  is said to satisfy the $\Delta _2$-condition
near infinity -- briefly, $A \in \Delta _2$ near infinity -- if it is
finite-valued and there exist constants $C>2$ and $t_0\geq 0$ such
that
\begin{equation}\label{delta2}
A(2t) \leq C A(t) \quad \hbox{for $t \geq t_0$.}
\end{equation}
Owing to equation \eqref{gen9},    condition \eqref{delta2} turns out to be equivalent to the existence of constants $C'>0$ and $t_1>0$ such that
\begin{equation}\label{delta2'}
a(2t) \leq C' a(t) \quad \hbox{for $t \geq t_1$.}
\end{equation}
 The function $A$ is said to satisfy the $\nabla _2$-condition
near infinity -- briefly, $A \in \nabla _2$ near infinity -- if there
exist constants $C>2$ and $t_0\geq 0$ such that
\begin{equation}\label{nabla2}
A(2t) > C A(t) \quad \hbox{for $t \geq t_0$.}
\end{equation}
One can show that
\begin{equation}\label{nabla2incr}
 \hbox{$A \in \nabla _2$ near infinity if and only if the function \,\, $t \mapsto \frac{A(t)}{t^{1+\varepsilon}}$ \,\,
is increasing for  $t\geq t_0$,}
\end{equation}
 for some  constants
$\varepsilon>0$ and $t_0\geq 0$.
Let us also note that
\begin{equation}\label{deltanabla}
\hbox{$A\in \Delta _2$ \, near infinity if and only if \, $\widetilde A \in \nabla _2$ \, near infinity.}
\end{equation}
A  Young function $A$ is said to
dominate another Young function $B$ near infinity   if there exist
constants $C>0$ and $t_0\geq 0$  such that
\begin{equation}\label{B.5bis}
B(t)\leq A(C t) \qquad \textrm{for \,\,\,} t\geq t_0\,.
\end{equation}
The functions  $A$  and  $B$ are called equivalent near infinity if
they dominate each other near infinity.
 \\ If any of the above definitions is satisfied with $t_0=0$, then it is
said to hold globally, instead of just near infinity.

\smallskip
\par
\par
Now, let $E$ be a measurable subset of $\rn$. 
We denote by $\mathcal{M}(E)$  the space of real-valued measurable functions on  $E$. The notation  $\mathcal{M}_+(E)$ is adopted for the subset of nonnegative functions in $\mathcal{M}(E)$. Similarly, the subscript $+$ attached to the notation of other spaces of real-valued functions will be used to denote the subset of nonnegative functions from those spaces.
\\ The Orlicz space $L^A (E)$ built upon a
Young function $A$  is the Banach function space of those functions 
$u\in \mathcal{M}(E)$ for which the
 Luxemburg norm
\begin{equation}\label{lux}
 \|u\|_{L^A(E)}= \inf \left\{ \lambda >0 :  \int_{\mathcal R }A
\left( \frac{|u |}{\lambda} \right) dx  \leq 1 \right\}\,
\end{equation}
is finite. In particular, $L^A (E)= L^p
(E)$ if $A(t)= t^p$ for some $p \in [1,
\infty )$, and $L^A (E)= L^\infty
(E)$ if $A(t)= 0$ for $t \in [0,1)$ and $A(t)= \infty$ for $t \in [1,\infty)$. 
\\ Moreover, we  denote by $L^A_{\rm loc}(E)$ the set of those functions from $\mathcal{M}(E)$ that belong to $L^A(F)$ for every bounded set $F\subset E$. 
  \\
 The H\"older type inequality
\begin{equation}\label{holder}
\|v\|_{L^{\widetilde A}(E)} \leq \sup_{u \in
L^A(E)}\frac{\displaystyle \int _E |u \,v
|\,dx}{\|u\|_{L^A(E)}}   \leq 2
\|v\|_{L^{\widetilde A}(E)}
\end{equation}
holds for every $u \in L^A(E)$ and $v\in
L^{\widetilde A}(E)$. 
\\ Denote by $|E|$ the Lebesgue measure of $E$, and assume that $|E|<\infty$.  Then
\begin{equation}\label{normineq}
L^{A}(E) \to L^{B}(E),
%
\end{equation}
if and only if the Young function $A$ dominates the Young function $B$ near
infinity. Here, and in what follows, 
where the arrow $\lq\lq \to"$ stands for continuous embedding. In
particular, 
\begin{equation}\label{normeq}
\hbox{$L^{A}(E) = L^{B}(E)$ (up to equivalent norms), if and only if $A$ and $B$
are equivalent near infinity.}
\end{equation}
\par
Given an open set $\Omega \subset \rn$ and a Young function $A$, 
the Orlicz-Sobolev space $W^{1,A}(\Omega)$ is defined as
$$W^{1,A}(\Omega) = \{u\in L^A(\Omega): \hbox{$u$ is weakly differentiable, and $|\nabla u|\in L^A(\Omega)$} \}.$$
The space $W^{1,A}(\Omega)$, equipped with the norm
\begin{equation}\label{OS}
 \|u\|_{W^{1,A}(\Omega)}= \|u\|_{L^A(\Omega)}\,+\,\|\nabla u\|_{L^A(\Omega)}
\end{equation}
is a Banach space. The space of those functions $u \in \mathcal M(\Omega)$ such that $u \in W^{1,A}(\Omega')$ for every bounded open set $\Omega ' \subset \subset \Omega$ will be denoted by $W^{1,A}_{\rm loc}(\Omega)$.
\par
We refer the reader to the monographs \cite{KR, RR1, RR2} for a comprehensive treatment of the topics of this section.

\section{Main results}\label{main}

We begin our analysis with a condition on  Young functions $A$ ensuring the local boundedness, as well as the continuity and  differentiability almost everywhere, of weakly monotone functions from the Orlicz-Sobolev space $W^{1,A}_{\rm loc}(\Omega)$, where $\Omega$ is an open subset of $\rn$, with $n \geq 2$. Its formulation involves the Young function $A_{n-1}$ associated with $A$ and $n$ as 
\begin{equation}\label{Apoincarebis}
A_{n-1}(t) = \begin{cases} A(t) & \quad \hbox{if $n=2$\,,} \\
\\
\displaystyle \left(t^{\frac {n-1}{n-2}} \int _t^\infty
\frac{\widetilde A(r)}{r^{1+\frac {n-1}{n-2}}}\,
dr\right)^{\widetilde {\quad}} & \quad \hbox{if $n\geq 3$\,,}
\end{cases}
\end{equation}
for $t\ge 0$. The integral on the right-hand side of \eqref{Apoincarebis} is convergent if and only 
the function $A$ fulfills condition \eqref{equiv2bis} -- see e.g. \cite[Lemma 4.1]{Cianchiibero}. This condition 
%
%
 will always come into play when dealing with the function $A_{n-1}$ for $n \geq 3$.
 \\ As mentioned in Section \ref{intro}, the function $A_{n-1}$ arises in a Poincar\'e type inequality for the oscillation of functions from  Orlicz-Sobolev spaces on the sphere, and  has a crucial role in the results to be presented.  Let us notice that the function $A_{n-1}(t)$ always dominates $A(t)$, and, if $n \geq 3$, it is equivalent to $A(t)$ whenever the latter grows faster than the  function $t^{n-1}$ in a suitable sense -- see equation \eqref{indexcond} below.
\\ In what follows, the notation $\medint _{B_r(x)} \cdots\, dz$ stands for $\frac 1{|B_r(x)|}\int _E \cdots\, dz$, where $B_r(x)$ is the ball, centered at $x\in \rn$, with radius $r>0$.

\begin{theorem}\label{orliczbound}
Let $A$ be a Young function.  Assume that either $n=2$, or $n \geq 3$ and $A$ fulfills  condition \eqref{equiv2bis}. Let $A_{n-1}$ be the Young function defined by \eqref{Apoincarebis}.
Let $u \in W_{\rm loc}^{1,A}(\Omega )$ be  a weakly monotone function.
Then $u \in L^\infty _{\rm loc}(\Omega)$, and there exists a constant $c=c(n)$ such that
\begin{equation}\label{bound1}
{\rm ess\,  osc}_{B_r(x)}u \leq c\, r\, A_{n-1}^{-1} \bigg( \medint _{B_{2r}(x)} A(|\nabla u|)\,dz \bigg)
\end{equation}
whenever $B_{2r}(x) \subset \subset \Omega$. Moreover, there exists a representative of $u$ that is   differentiable  a.e. in $\Omega$. 
\end{theorem}

\begin{remark}\label{onninen} {\rm The a.e. differentiability of weakly differentiable functions $u \in W_{\rm loc}^{1,A}(\Omega )$ under assumption \eqref{equiv2bis} can also be derived from  \cite[Theorem 1.2]{O}, \color{black}
 via an inclusion relation between Orlicz and Lorentz spaces established in \cite{KKM}.  Here, we present a self-contained proof, that just relies upon Orlicz spaces techniques, and contains some preliminary steps of  use for our subsequent results.}
\end{remark}

\smallskip

More precise information about the set of points of continuity of any weakly differentiable function $u \in W_{\rm loc}^{1,A}(\Omega )$ can in fact be provided under assumption \eqref{equiv2bis}. It turns out that any such function has a representative whose restriction to the complement in $\Omega$ of an exceptional set of (suitably defined) vanishing capacity is continuous. This is the content of Theorem \ref{orliczcapacity} below. 
\\ The relevant capacity generalizes the standard $C_{p, 1}$ capacity associated with the Sobolev space $W^{1,p}_{\rm loc}(\Omega)$, and depends on the Young function $A$ and on the dimension $n$ of the ambient space $\rn$ of $\Omega$. It can be defined as follows. 
\\
Let $\Psi : [0, \infty) \to [0, \infty)$ be a continuous function.
%
%
 Consider the Riesz type operator defined as
\begin{equation}\label{IPsi}
I_\Psi f (x) = \int _{\rn}\frac{f(y)}{|x-y|^{n} \Psi(1/|x-y|)}\, dy
\quad \hbox{for $x \in \rn$}
\end{equation}
for $f \mathcal \in \mathcal M_+(\rn)$.
 The associated  capacity
$C_{\Psi , 1}$ of a set $E \subset \rn$ is given by
\begin{equation}\label{cap}
C_{\Psi , 1} (E) = \inf\bigg\{\int _{\rn} f(x)\, dx: \, f\in \mathcal
M_+(\rn), \, I_\Psi f (x)\geq 1 \,\,\hbox{for $x \in E$}\bigg\}.
\end{equation}
Note that, with the choice $\Psi(t) = t^\alpha$, where $\alpha \in
(0, n)$, the operator $I_\Psi$ reproduces the classical Riesz potential $I_\alpha$ given by
\begin{equation}\label{riesz}
I_{\alpha} f (x) = \int _{\rn}\frac{f(y)}{|x-y|^{n-\alpha}}\, dy
\quad \hbox{for $x \in \rn$.}
\end{equation}
 Hence, 
the capacity $C_{\Psi , 1}$ agrees with the standard $C_{\alpha ,
1}$ capacity associated with the operator $I_{\alpha}$
\cite{AH}.

\begin{theorem}\label{orliczcapacity}
%
%
Let $A$ be a Young function.  Assume that either $n=2$, or $n \geq 3$ and $A$ fulfills  condition \eqref{equiv2bis}. Let $A_{n-1}$ be the Young function defined by \eqref{Apoincarebis}. Assume that $\sigma : [0, \infty) \to
[0, \infty)$ is a continuous function such that
\begin{equation}\label{ott21}
\int ^\infty \frac{A_{n-1}(\lambda t)}{t\sigma (t) A_{n-1}(t)}\, dt = \infty \qquad \hbox{for every $\lambda >0$.}
\end{equation}
Let $\Psi : [0, \infty ) \to [0, \infty )$ be the function given by
\begin{equation}\label{Psi}
\Psi (t) = \sigma (t) A_{n-1}(t) \quad \hbox{for $t \geq0$.}
\end{equation}
Then   every weakly monotone function $u \in W_{\rm loc}^{1,A}(\Omega )$
admits a representative  whose restriction to the complement in $\o$ of an exceptional set  of $C_{\Psi,1}$-capacity zero is continuous.
\end{theorem}

Assumption \eqref{ott21} takes a simpler form in the special case when the function $A \in \Delta _2$ near infinity. Indeed,  $A_{n-1}\in \Delta _2$ near infinity as well in this case -- see Proposition \ref{delta2lem}, Section \ref{secproofs}. Hence, $A_{n-1}(\lambda t)$ and  $A_{n-1}(t)$ are bounded by each other, up to positive multiplicative constants depending on $A$, $n$ and $\lambda$, for large $t$. This is stated in the following corollary.

\begin{corollary}\label{orliczcapacitycor}
%
%
Let $A$ be a Young function such that  $A \in \Delta_2$ near infinity. Assume that either $n=2$, or $n \geq 3$ and $A$ fulfills  condition \eqref{equiv2bis}. Let $A_{n-1}$ be the Young function defined by \eqref{Apoincarebis}. Assume that $\sigma : [0, \infty) \to
[0, \infty)$ is a continuous function such that
\begin{equation}\label{dini}
\int ^\infty \frac{dt}{t \sigma (t)} = \infty.
%
\end{equation}
Let $\Psi$ be the function defined as in \eqref{Psi}.
Then   every weakly monotone function $u \in W_{\rm loc}^{1,A}(\Omega )$
admits a representative  whose restriction to the complement in $\o$ of an exceptional set  of $C_{\Psi,1}$-capacity zero is continuous.
\\ In particular, condition \eqref{dini} holds with $\sigma (t)=1$, and hence $u$ enjoys this property outside a set of $C_{A_{n-1},1}$-capacity zero.
\end{corollary}

\begin{remark}\label{improveManfredi}  {\rm Corollary \ref{orliczcapacitycor} not only recovers, but  somewhat allows for improvements of  the result \eqref{p>n-1} of \cite{M} in the case when $A(t)=t^p$ for some $p\in (n-1,n)$. Indeed,  it tells us that any weakly monotone function in $W^{1,p}_{\rm loc}(\Omega)$, with $p\in (n-1, n)$, admits a representative that is continuous  outside a subset of $\Omega$ having zero  $C_{\sigma (t)t^p, 1}$-capacity, for any function $\sigma$ fulfilling \eqref{dini}. Possible choices are thus, for instance, $\sigma (t) = \log (1+ t)$, $\sigma (t)= \log (1+t) \log (1+ \log (1+t))$, etc.. 
\\ Let us also notice that, since no restriction is imposed on $A$ if $n=2$,  in this case Theorems \ref{orliczbound}  and \ref{orliczcapacity}   also hold  if $A(t)=t$, namely when $W^{1,A}_{\rm loc}(\Omega)= W^{1,1}_{\rm loc}(\Omega)$. In particular, this shows that the    result \eqref{p>n-1} of \cite{M} is still valid for the endpoint value $p=n-1=1$ when $n=2$.
}
\end{remark}

Capacities can be dismissed in the description of the size of the exceptional set   of possible discontinuity points of a weakly monotone function in $W_{\rm loc}^{1,A}(\Omega )$, with $A$ satisfying \eqref{equiv2bis}. Indeed,
its size   can  be estimated in terms of  Hausdorff measures. Given a continuous, increasing  function $h : [0, \infty) \to [0,
\infty)$ such that $h(0)=0$, the classical $h$-Hausdorff measure
$\mathcal H ^{h(\cdot)} (E)$ of a set $E \subset \rn$ is defined as
\begin{equation}\label{hausdorff}
\mathcal H ^{h(\cdot)} (E) = \lim _{\varepsilon \to 0^+} \inf \bigg\{\sum
_{j=1}^\infty h ({\rm diam}(K_j)): E \subset \cup _{j=1}^\infty
K_j\,,{\rm diam}(K_j) \leq \varepsilon\bigg\}.
\end{equation}
If $h(t) = t^\beta$ for some some $\beta >0$, then $\mathcal H ^{h(\cdot)}$
agrees (up to a multiplicative constant) with the standard
$\beta$-dimensional Hausdorff measure $\mathcal H ^\beta$.
\\ 
Note that we may assume, without loss of generality, that 
%
\begin{equation}\label{immaterial}r \mapsto \frac{h(r)}{r^n}
 \quad \hbox{is a non-increasing function.}
\end{equation}
%
Indeed, \cite[Proposition 5.1.8]{AH} tells us what follows. If $\liminf _{r\to 0^+}\tfrac{h(r)}{r^n} =0$, then $\mathcal H^{h(\cdot)}(E)=0$ for every set $E\subset \mathbb R^n$. If $\liminf _{r\to 0^+}\tfrac{h(r)}{r^n}>0$, then there exists another continuous increasing function $\overline h$ such that \eqref{immaterial} is satisfied with $h$ replaced by $\overline h$, and moreover $\mathcal H^{h(\cdot)} (E)$ and $\mathcal H^{\overline {h}(\cdot)} (E)$ are bounded by each other for every set $E\subset \mathbb  R^n$, up to multiplicative constants independent of $E$.

%

\smallskip
\par
Our estimate of the exceptional set of weakly monotone Orlicz-Sobolev functions via Hausdorff measures   involves a function $\Psi$ satisfying the assumptions of Theorem \ref{orliczcapacity}, and such that, in addition, 
\begin{equation}\label{febbraio1}
t \mapsto t^n\Psi (1/t) \quad \hbox{is a
non-decreasing function decaying to $0$ as $t\to 0^+$.}
\end{equation}
Obviously, this assumption ensures that the derivative $d\big(\frac {-1}{s^n\Psi (1/s)}\big)$ defines  a positive measure on $(0, \infty)$.

\begin{theorem}\label{continuityhausdorff}
Let $A$ be a Young function.  Assume that either $n=2$, or $n \geq 3$ and $A$ fulfills  condition \eqref{equiv2bis}.  Assume that the function  $h$ is as above. Let   $\Psi$ be as
in Theorem \ref{orliczcapacity}.  Assume, in addition,  that $\Psi$ satisfies condition \eqref{febbraio1}. If
\begin{equation}\label{hB}
\int _0 h(s)\, d\Big(\frac {-1}{s^n\Psi (1/s)}\Big) \, < \infty,
\end{equation}
then   every weakly monotone function $u \in W_{\rm loc}^{1,A}(\Omega )$
admits a representative  whose restriction to the complement in $\o$ of an exceptional set  of $\mathcal H ^{h(\cdot)}$-measure zero is continuous.
\end{theorem}

We are now ready to state the main result of this paper, concerning  the everywhere continuity of Orlicz-Sobolev weakly monotone functions. It asserts that, 
if assumption \eqref{equiv2bis}  is properly strengthened, then any weakly monotone function from $W^{1,A}_{\rm loc}(\Omega)$ has a representative which is  continuous in the whole of $\Omega$. The assumption to be imposed is \eqref{ours}, 
%
%
with $A_{n-1}$  defined by \eqref{Apoincarebis}. In fact, under  assumption \eqref{ours} any weakly monotone function  from $W^{1,A}_{\rm loc}(\Omega)$ is
locally uniformly continuous with a modulus of continuity depending only on $A$ and $n$. This  modulus of continuity 
 $\omega : [0, \infty) \to [0, \infty)$ is defined as
\begin{equation}\label{modulus}
\omega (r) = r B^{-1}(r^{-n}) \quad \hbox{for $r
>0$,}
\end{equation}
where
%
%
%
 $B: [0, \infty) \to [0, \infty)$ is the function given by
\begin{equation}\label{G} B(t)= t^n \int _0^t \frac {A_{n-1}(s)}{s^{1+n}}\, ds
\quad \hbox{for $t >0$.}
\end{equation}
The space of functions in $\o$ that are locally uniformly continuous with modulus of continuity not exceeding $\omega$ will be denoted by $C^{\omega (\cdot )}_{\rm loc}(\Omega)$.
\\
Observe that, since we are dealing with  properties of functions from the local
Orlicz-Sobolev spaces $W^{1,A}_{\rm loc}(\Omega)$, the function $A$ can be modified, if
necessary,  near $0$ in such a way that
\begin{equation}\label{febbraio4}
\int _0\frac {A_{n-1}(t)}{t^{1+n}}\, dt < \infty.
\end{equation}
Owing to property \eqref{normeq}, a modification of this kind leaves the space $W^{1,A}_{\rm loc}(\Omega)$ unchanged. With condition \eqref{febbraio4} in force, the function $B$ is well defined. Also, it can be verified that it is  a Young function, whence its inverse $B^{-1}$ is well defined as well.

\begin{theorem}\label{febbraio0}
Let $n \geq 2$, and let $A$ be a Young function fulfilling condition \eqref{ours}, 
%
%
with $A_{n-1}$  defined by \eqref{Apoincarebis}.
Then   every weakly monotone function $u \in W_{\rm loc}^{1,A}(\Omega )$
admits a continuous representative. Moreover,
$u \in C^{\omega (\cdot )}_{\rm loc}(\Omega)$, where $\omega$ is given by \eqref{modulus}. 
\end{theorem}

\begin{remark}\label{remmodulus}
{\rm If $n=2$, 
assumption \eqref{ours} agrees with
\begin{equation}\label{mainhp2}
\int ^\infty \frac{A(t)}{t^3}\, dt =
\infty \,,
\end{equation}
namely with \eqref{Iwcond}, since $A_1=A$. On the other hand, if $n \geq 3$, assumption \eqref{ours} 
is equivalent to 
\begin{equation}\label{mainhp}
 \int ^\infty t^{\frac{1-n}{n-2}}
\bigg(\int _t^\infty \frac{\widetilde A (s)}{s^{1+\frac
{n-1}{n-2}}}ds\bigg)^{1-n}\, dt  = \infty \quad \hbox{if $n\geq 3$\,.}
\end{equation}
Indeed, by \cite[Lemma 4.1]{Cianchiibero}, condition \eqref{ours}
is equivalent to
$$\int ^\infty \bigg(\frac s{\widetilde {A_{n-1}}(s)}\bigg)^{n-1}ds
= \infty,$$ and, in view of  definition \eqref{equiv2bis}, the latter coincides with \eqref{mainhp}.}
\end{remark}

Under the additional assumption that $A \in \Delta_2$ near infinity, condition \eqref{ours} can be reformulated, for $n \geq 3$, in a form which only involves $A$, and avoids explicit reference to $A_{n-1}$. This fact is a consequence of \cite[Lemma 3.3]{CarozzaCianchi}, and is enucleated in the next result.

\begin{corollary}\label{maindelta2}
Let $n \geq 3$, and let $A$ be a Young function such that $A \in \Delta_2$ near infinity, and
\begin{equation}\label{condelta2}
\int ^\infty \Big(\frac t{A(t)}\Big)^{\frac 2{n-2}} \bigg(\int _t^\infty \Big(\frac s{A(s)}\Big)^{\frac 1{n-2}}\,ds \bigg)^{-n}\, dt = \infty\,.
\end{equation}
Then   every weakly monotone function $u \in W_{\rm loc}^{1,A}(\Omega )$
admits a continuous representative. Moreover,
$u \in C^{\omega (\cdot )}_{\rm loc}(\Omega)$, where $\omega$ is given by \eqref{modulus}.
\end{corollary}

Let us next mention a standard situation when assumption \eqref{ours} reduces to \eqref{Iwcond}, even for $n\geq 3$. Suppose, for simplicity, that $A$ is finite-valued. This is of course the only nontrivial case, since $W^{1,A}_{\rm loc}(\Omega)= W^{1,\infty}_{\rm loc}(\Omega)$ if $A$ jumps to infinity, and every function in the latter space is locally Lipschitz continuous. When $n\geq 3$, 
the function $A_{n-1}$  is equivalent to $A$ if and only if its lower Boyd index at infinity, defined as 
\begin{equation}\label{index}
i(A)= \lim _{\lambda \to \infty} \frac{\log \Big(\liminf _{t\to
\infty }\tfrac{A(\lambda t)}{A(t)}\Big)}{\log \lambda },
\end{equation} 
satisfies 
\begin{equation} \label{indexcond}
i(A) >n-1\,,
\end{equation}
see \cite[Lemma 2.3]{Strom} and, for more details, \cite[Proposition 4.1]{CianchiMusil}. Thus, under assumption \eqref{indexcond}, condition \eqref{ours} agrees with \eqref{Iwcond}. As a consequence, the following corollary of Theorem \ref{febbraio0} holds. Note that, since assumption \eqref{incr} implies \eqref{indexcond}, this corollary recovers, in particular, the results of \cite{IM} and \cite{KKMOZ} mentioned in Section \ref{intro}.

\begin{corollary}\label{recover}
Let $A$ be a Young function fulfilling condition \eqref{Iwcond}. Assume that either $n=2$, or $n \geq 3$ and \eqref{indexcond} holds. Then every weakly monotone function in 
$W_{\rm loc}^{1,A}(\Omega )$ is continuous in $\Omega $.
\end{corollary}

\color{black}
\color{black}

\section{Proofs of the main results}\label{secproofs}

The proof of Theorem \ref{orliczbound} requires  a preliminary result on smooth approximation of functions from an Orlicz-Sobolev  space $W^{1,A}(\Omega)$. Such an approximation is not possible in  norm, unless $A$ satisfies a $\Delta _2$-condition. However, standard convolution with a sequence of smooth kernels always provides us with an approximating sequence whose Dirichlet integrals associated with $A$ converge to the corresponding  Dirichlet integral of the limit function.

\begin{lemma}\label{approx}
Let $A$ be a Young function, and let $u  \in W^{1,1}(\Omega)$ be
such that $\int _\Omega A(|\nabla u|)\,dx < \infty$. Let $\{u_k\}$ be
a sequence of convolutions of $u$ with mollifiers $\varrho _k$, namely
\begin{equation}\label{moll} u_k = u * \varrho  _k \,,\end{equation}
where $\varrho _k \in C^\infty_0(\rn)$, ${\rm supp} \varrho _k \subset B_{1/k}(0)$, $\varrho  _k \geq 0$ and $\int _{\rn}\varrho _k\, dx=1$ for $k \in \mathbb N$.
Then
(up to subsequences)
\begin{equation}\label{approx1}
\lim _{k \to \infty} u_k  = u   \quad \hbox{at every Lebesgue point  of $u$,}
\end{equation}
\begin{equation}\label{approx1grad}
\lim _{k \to \infty}  \nabla u_k =  \nabla u  \quad \hbox{at every Lebesgue point  of $\nabla u$,}
\end{equation}
and
\begin{equation}\label{approx2}
\lim _{k \to \infty}  \int _E A(|\nabla u_k|)\,dx  = \int _E A(|\nabla u|)\, dx   \quad
 \hbox{for every measurable set $E \subset \Omega$.}
\end{equation}
\end{lemma}
\par\noindent
{\bf Proof}. Properties  \eqref{approx1} an \eqref{approx1grad} are classical. As far as
\eqref{approx2} is concerned, given any measurable set $F \subset E$, the Hardy-Littlewood inequality ensures that 
\begin{align}\label{approx4}
    \int _F A(|\nabla u_k|)\,dx  \leq \int _0^{|F|}A(|\nabla u_k|^*(r))\,
    dr \,,
\end{align}
where the asterisque $\lq\lq *"$ stands for decreasing rearrangement. Moreover, 
a rearrangement inequality for convolutions \cite{oneil},
      \begin{align}\label{approx3}
     \int _0^\tau |\nabla u_k|^{*} (s)\, ds & \leq \int_0^\tau|\nabla u|^{*}(s)\, ds \int _0^\tau \varrho _k^*(s)\, ds
      + \tau\, \int _\tau ^\infty |\nabla u|^*(s) \varrho _k^*(s)\, ds
    \\ \nonumber & \leq \int_0^\tau|\nabla u|^{*}(s)\, ds \int _0^\tau \varrho _k^*(s)\, ds + \tau \,|\nabla u|^*(\tau ) \int _\tau ^\infty  \varrho _k^*(s)\, ds
   \\ \nonumber & \leq \int_0^\tau|\nabla u|^{*}(s)\, ds \int _0^\infty \varrho _k^*(s)\, ds
    = \int_0^\tau|\nabla u|^{*}(s)\, ds \int _{\rn} \varrho _k (x)\, dx  \\ \nonumber &  =
    \int_0^\tau|\nabla u|^{*}(s)\, ds \quad \quad \hbox{for $\tau \geq 0.$}
      \end{align}
    Note that, in the last but one inequality, we have made use of the fact that  $\int _0^\infty \varrho _k^*(s) \, ds = \int _{\rn } \varrho _k (x)\, dx =1$ for every $k \in \N$. Inequality \eqref{approx3}, via \cite[Proposition 2.1]{ALT}, tell us that
    \begin{align}\label{nov7}
    \int _0^{|F|} A(|\nabla u_k|^*(s))\,
    ds \leq \int _0^{|F|} A(|\nabla u|^*(s))\,ds\,.
    \end{align}
    Coupling inequalities \eqref{approx4} and \eqref{nov7} yields
      \begin{align}\label{nov8}
    \int _F A(|\nabla u_k|)dx  \leq  \int _0^{|F|} A(|\nabla u|^*(s))\,ds\,.
    \end{align}
    Inequality \eqref{nov8} entails that the sequence $\{A(|\nabla u_k|)\}$ is equi-integrable over $E$, since
    $$\int _0^{|F|} A(|\nabla u|^*(s))\,ds \leq \int _0^{|\Omega|} A(|\nabla u|^*(s))\,ds \leq \int _\o A(|\nabla u|)\, dx\,.$$
    Hence, \eqref{approx2} follows via \eqref{approx1grad} and Vitali's
    convergence theorem. \qed
    
    \medskip
\par\noindent
We are now ready to prove Theorem \ref{orliczbound}.

\medskip
\par\noindent
{\bf Proof of Theorem \ref{orliczbound}}.
On replacing, if necessary, $\Omega$ with a bounded open subset, we may  suppose, without loss of generality, that $u \in W^{1,A}(\Omega)$. Moreover, we assume, for the time being, that 
\begin{equation}\label{finiteint}
\int _{\Omega}A(|\nabla u|)\,dz < \infty.
\end{equation}
Let $u_k$ be the sequence appearing in Lemma \ref{approx}. Given any 
$x_0 \in \Omega$,  let  $R>0$ be such that $B_R(x_0) \subset \subset
\Omega$, and let $r \in (0, R)$. By  \cite[Lemma 7.4.1]{IM}, for every $\delta
>0$ and any Lebesgue points  $x,y \in B_r(x_0)$ of $u$, there exists $\overline k = \overline k (x, y, \delta , r, R)$ such that
\begin{equation}\label{osc}
|u_k(x)\,-\,u_k(y)|\leq 2 \delta + {\rm osc} _{S_\tau (x_0)} u_k
\end{equation}
if $k \geq \overline k $ and $\tau \in [r, R]$. Here, $S_\tau (x_0)$ denotes the $(n-1)$-dimensional sphere in $\rn$ centered at $x_0$,  with radius $\tau$, and $$ {\rm osc} _{S_\tau (x_0)} u_k = {\rm sup} _ {S_\tau (x_0)} u_k - {\rm inf}_ {S_\tau (x_0)} u_k.$$
An Orlicz-Sobolev
Poincar\'e type inequality on the $(n-1)$-dimensional sphere $S_\tau (x_0)$
\cite[Theorem 4.1]{CarozzaCianchi} (see also \cite{AlCi, Cianchi_ASNS, Mazya, Talenti} for related results) tells us that, if either $n=2$, or $n \geq 3$ and \eqref{equiv2bis} holds, then
\begin{equation}\label{orlicz3}
{\rm osc}_{S_\tau (x_0)}u_k \leq C \tau A_{n-1}^{-1} \bigg( \tau ^{1-n} \int
_{S_\tau (x_0)} A(|\nabla u_k|) d\mathcal H ^{n-1} \bigg),
\end{equation}
for some constant $C=C(n)$, and for $\tau >0$.
Thanks to \eqref{osc} and \eqref{orlicz3}, 
\begin{equation}\label{orlicz17}
\frac{1}{C\tau}|u_k(x)\,-\,u_k(y)|\leq \frac{2 \delta}{C\tau} +
A_{n-1}^{-1} \bigg( \tau^{1-n} \int _{S_\tau (x_0)} A(|\nabla u_k|)
d\mathcal H ^{n-1} \bigg)\,,
\end{equation}
if $k \geq \overline k$ and $\tau\in [r, R]$.
 Given $\alpha \in (0,1)$, inequality \eqref{orlicz17}
can be rewritten as
\begin{equation}\label{orlicz5'}
\frac{\alpha}{C\tau}|u_k(x)\,-\,u_k(y)|\leq \frac{2\alpha
\delta}{C\tau (1-\alpha)}(1-\alpha) +\alpha A_{n-1}^{-1}
\bigg(\tau^{1-n} \int _{S_\tau (x_0)} A(|\nabla u_k|) d\mathcal H
^{n-1} \bigg)\,.
\end{equation}
Hence, by the convexity of the 
 function $A_{n-1}$,  
\begin{equation}\label{orlicz5''}
A_{n-1}\bigg(\frac{\alpha}{C\tau}|u_k(x)\,-\,u_k(y)|\bigg)\leq
(1-\alpha)A_{n-1}\bigg(\frac{2\alpha \delta}{C\tau (1-\alpha)}\bigg)
+\alpha  \tau^{1-n} \int _{S_\tau (x_0)} A(|\nabla u_k|) \,d\mathcal H
^{n-1} \,.
\end{equation}
Now, fix any Lebesgue point $t\in (r,\,R)$ for the function $t\mapsto \int
_{S_t(x_0)} A(|\nabla u|) \,d\mathcal H ^{n-1}$. Note that this
function belongs to $L^1(0,R)$, since  $\int _{B_R(x_0)} A(|\nabla
u|) \,dx <\infty$. Given any  number $\varepsilon
>0$ such that $(t-\ep , t+\ep) \subset (r,R)$, multiply through by $\tau ^{n-1}$ inequality \eqref{orlicz5''},  and integrate over $(t-\varepsilon,
t+\varepsilon)$ to obtain
\begin{multline}\label{orlicz6'}
\int_{t-\varepsilon}^{t+\varepsilon}
\tau^{n-1}A_{n-1}\bigg(\frac{\alpha}{C\tau}|u_k(x)\,-\,u_k(y)|\bigg)d\tau
\\ \leq 2\varepsilon (1-\alpha)R^{n-1}A_{n-1}\bigg(\frac{2\alpha
\delta}{Cr (1-\alpha)}\bigg)+ \alpha
\int_{t-\varepsilon}^{t+\varepsilon} \int _{S_\tau(x_0)} A(|\nabla
u|)\,d\mathcal H ^{n-1}\, d\tau\,.
%
\end{multline}
Passing to the limit  as $k\to \infty$ in inequality \eqref{orlicz6'}, making use of equations
\eqref{approx1} and \eqref{approx2}, and then passing to the limit as
$\delta \to 0$ yield
\begin{equation}\label{orlicz7}
\int_{t-\varepsilon}^{t+\varepsilon}\frac{\tau^{n-1}}{\alpha
}A_{n-1}\bigg(\frac {{\alpha |u(x)\,-\,u(y)|}}{C\tau}\bigg)d\tau
\leq \int_{t-\varepsilon}^{t+\varepsilon}  \int _{S_\tau(x_0)}
A(|\nabla u|)\,d\mathcal H ^{n-1}\,d\tau \;.
\end{equation}
Since $A_{n-1}$ is a Young function, the function
$A_{n-1}(\alpha )/\alpha$ is increasing in $\alpha$. One can then pass to the limit as $\alpha \to 1^-$, and make use  of the
monotone convergence theorem in the integral on the left-hand side
of \eqref{orlicz7} to deduce that inequality \eqref{orlicz7} continues to hold for $\alpha =1$. On dividing through by $2\varepsilon $  the
resulting inequality, and letting $\varepsilon\to 0$ we conclude
that
\begin{equation}\label{orlicz4}
t^{n-1}A_{n-1}\bigg(\frac {|u(x) - u(y)|}{Ct}\bigg) \leq \int
_{S_t(x_0)} A(|\nabla u|)\, d\mathcal H ^{n-1}
\end{equation} 
 for all Lebesgue points $x,y
\in B_r(x_0)$ of $u$,  and for a.e. $t \in [r, R]$.
Next, observe that
the function
\begin{equation}\label{noninc} 
t \mapsto  t^{n-1}
A_{n-1}(1/t) \quad \hbox{is decreasing.}
\end{equation}
  Indeed, property \eqref{noninc} is equivalent to the fact that the function 
  \begin{equation}\label{noninc2} 
  t\mapsto \frac{A_{n-1}(t)}{t^{n-1}} \quad \hbox{is increasing.}
\end{equation}
  If $n=2$, then  property \eqref{noninc2} just holds because $A_{n-1}$ is a Young function. If $n \geq 3$, then by \eqref{equivconj} property \eqref{noninc2} is in turn equivalent to the fact that the function 
  \begin{equation}\label{noninc3} 
  t\mapsto \frac{\widetilde {A_{n-1}}(t)}{t^{(n-1)'}} \quad \hbox{is decreasing.}
\end{equation}
  Property \eqref{noninc3} trivially holds, since 
  $$\frac{\widetilde {A_{n-1}}(t)}{t^{(n-1)'}} = 
  \int _t ^\infty
\frac{\widetilde A (r)}{r^{1 + \frac {n-1}{n-2}}} \, dr \quad \hbox{for $t>0$.}
$$
%
 Now, assume that $2r < R$. On integrating \eqref{orlicz4} over $(r,
2r)$ and making use of  \eqref{noninc}
we deduce that
\begin{align*}
r (2r)^{n-1} A_{n-1}\bigg(\frac {|u(x)-u(y)|}{2Cr}\bigg) & \leq
\int _r ^{2r} t^{n-1}A_{n-1}\bigg(\frac {|u(x) -
u(y)|}{Ct}\bigg)\, dt
 \leq \int _r ^{2r} \int _{S_t(x_0)} A(|\nabla
u|) \,d\mathcal H ^{n-1}\, dt
 \\
\nonumber & \leq \int
_0 ^{2r} \int _{S_t(x_0)} A(|\nabla u|)\, d\mathcal H ^{n-1}\, dt
 = \int _{B_{2r}(x_0)} A(|\nabla u|) \,dz < \infty\,,
\end{align*}
whence \eqref{bound1} follows. Note that here we have also made use of  property \eqref{kt-1} with $A$ replaced by  $A_{n-1}$.
\\ Next,  let $\widehat u : \Omega \to (-\infty , \infty]$ be the function defined by
\begin{equation}\label{utilde}
\widehat u (x) = \limsup_{r \to 0^+} \medint _{B_r(x)} u(z) \, dz \quad \hbox{for $x\in \Omega$. }
\end{equation}
By Lebesgue differentiation theorem,
\begin{equation}
u(x) = \widehat u(x) \quad \hbox{for a.e. $x\in \Omega$.}
\end{equation}
 Given any $x,y \in B_r(x_0)$, let $\rho, \sigma >0$ be such that $B_\rho(x) \subset B_r(x_0)$ and $B_\sigma (y) \subset B_r(x_0)$. Owing to inequality \eqref{bound1}, 
\begin{equation}\label{medie}
\bigg| \medint _{B_\rho(x)} u \, dz - \medint _{B_\sigma(y)} u \, dz\bigg| \leq C r A_{n-1}^{-1} \bigg( \medint _{B_{2r}(x_0)} A(|\nabla u|)\,dz \bigg)\,.
\end{equation}
Passing to the limit in \eqref{medie} first as $\rho \to 0^+$, and then as $\sigma \to 0^+$ tells us that 
\begin{equation}\label{osctilde}
| \widehat u(x) -  \widehat u(y)| \leq C r A_{n-1}^{-1} \bigg( \medint _{B_{2r}(x_0)} A(|\nabla u|)\,dz \bigg)\,.
\end{equation}
Inequality \eqref{osctilde} ensures that $ \widehat u$ is continuous at every Lebesgue point of the function  $A(|\nabla u|)$, and hence a.e. in $\Omega$. Moreover, an application of \eqref{osctilde} with $y=x_0$ and $r=|x-x_0|$ yields 
$$\frac{|  \widehat u(x) -  \widehat u(x_0)|}{|x -x_0|} \leq C  A_{n-1}^{-1} \bigg( \medint _{B_{2|x -x_0|}(x_0)} A(|\nabla u|)\,dz \bigg)\,.$$
 Thus,
$$\limsup _{x \to x_0}\frac{|  \widehat u(x) -  \widehat u(x_0)|}{|x -x_0|} < \infty$$
if $x_0$ is any Lebesgue point  of the function  $A(|\nabla u|)$. Hence, the a.e. differentiability of $\widehat u$ a.e. in $\Omega$ follows, via a classical result  by Stepanoff \cite{Step}.  
\\Finally, if the temporary condition \eqref{finiteint} does not hold, then, however, it does hold with $u$ replaced with  $\lambda u$ for a suitable $\lambda >0$. The above argument then applies to $\lambda u$, and hence the continuity and the a.e. differentiability of $u$ still follows. As for inequality \eqref{bound1}, it trivially continues to hold even if \eqref{finiteint} fails, since its right-hand side is infinite in this case. \color{black}
\qed

\medskip\par\noindent
The content of the next lemma is a basic property of the capacity defined as in \eqref{cap}, to be used in the proof of Theorem \ref{orliczcapacity}.

\begin{lemma}\label{caplemma}
Let $\Psi : [0, \infty) \to [0, \infty)$ be a continuous function, and let $C_{\Psi,1}$ be the capacity defined by \eqref{cap}.
%
If $f \in
L^1_+(\rn)$, then
\begin{equation}\label{cap0} C_{\Psi,1} \big(\{I_\Psi  f  = \infty\}\big)
=0.
\end{equation}
\end{lemma}
\par\noindent
{\bf Proof}. Given  $f \in L^1_+(\rn)$, we have that
\begin{align}\label{febbraio5}
C_{\Psi,1} \big(\{I_\Psi  f  = \infty\}\big)\leq  C_{\Psi,1}\big(\{I_\Psi
f>\lambda \}\big) = C_{\Psi,1}\big(\{I_\Psi (f/\lambda)>1\}\big) \leq
\frac 1\lambda \int _{\rn} f(x)\, dx,
\end{align}
for every $\lambda >0$. Hence, equation \eqref{cap0} follows on
letting $\lambda$ go to infinity. \qed

\medskip
\par\noindent
{\bf Proof of Theorem \ref{orliczcapacity}}. On replacing, if necessary, $u$ by $\lambda u$ for a suitable $\lambda >0$, we may assume, without loss of generality, that condition \eqref{finiteint} is in force.
Given $x_0\in \Omega$, let $R>0$ be such that $B_R(x_0) \subset \subset \Omega$. Set
$$ {\rm ess \,osc _{B_r(x_0)}} u = {\rm ess \,sup}_{B_r(x_0)} u - {\rm ess \,inf}_{B_r(x_0)} u$$
for $r \in (0, R)$.
 Inequality
\eqref{orlicz4} implies that 
\begin{equation}\label{gen1}
t^{n-1}A_{n-1}\bigg(\frac {{\rm ess\,osc}_{B_r(x_0)}u}{Ct}\bigg) \leq \int
_{S_t(x_0)} A(|\nabla u|) \,d\mathcal H ^{n-1}
\end{equation}
for a.e. $t \in [r, R]$.
Multiplying through equation \eqref{gen1}
by $\frac 1{t^n \sigma (1/t) A_{n-1}(1/t)}$, and  
integrating the resulting equation over $(r, R)$ yields
\begin{align}\label{feb1}
\int _r^R \frac{1}{t\sigma (1/t)}\frac{A_{n-1}({\rm ess\,osc}_{B_r(x_0)}u/(tC))}{A_{n-1}(1/t)}\, dt
& \leq \int _r^R \frac{1}{t^n \sigma (1/t)A_{n-1}(1/t)} \int
_{S_t(x_0)} A(|\nabla u|) \, d\mathcal H ^{n-1}\, dt
\\ \nonumber & \leq
 \int _0^R \frac{1}{t^n \sigma (1/t)A_{n-1}(1/t)} \int
_{S_t(x_0)} A(|\nabla u|)\, d\mathcal H ^{n-1}\, dt.
\end{align}
A change of variables to polar coordinates tells us that
\begin{align}\label{feb2}
\int _0^R &\frac{1}{t^n \sigma (1/t)A_{n-1}(1/t)} \int _{S_t(x_0)}
A(|\nabla u|) \,d\mathcal H ^{n-1}\, dt \\ \nonumber & = \int _{B_R(x_0)}\frac{
A(|\nabla u|)}{|x-x_0|^n \sigma (1/|x-x_0|)A_{n-1}(1/|x-x_0|)}\,dx
 = I_{\Psi}\big(A(|\nabla u|) \chi_{B_R(x_0)}\big)\,.
\end{align}
By Lemma \ref{caplemma}, there exists a set $E\subset \Omega$ such that $C_{\Psi,1}(E)=0$, and $I_{\Psi}\big(A(|\nabla u|)
\chi_{B_R(x_0)}\big)$ is finite for every $x_0\in \Omega\setminus E$.
 We claim that
 \begin{equation}\label{feb11}
 \lim _{r\to 0^+} {\rm ess \,osc}_{B_r(x_0)}u =0
 \end{equation}
 for any such $x_0$. Assume, by contradiction, that \eqref{feb11}
 fails, and hence ${\rm ess\,osc}_{B_r(x_0)}u\geq \lambda $  for some $\lambda >0$. 
Equations  \eqref{feb1} and \eqref{feb2}  imply that
\begin{align*}
C' \geq  \int _r^R  \frac{1}{t\sigma (1/t)}\frac{A_{n-1}({\rm
ess\,osc}_{B_r(x_0)}u/(tC))}{A_{n-1}(1/t)}\, dt & \geq \int _r^R
\frac{1}{t\sigma (1/t)}\frac{A_{n-1}(\lambda/(tC))}{A_{n-1}(1/t)}\,
dt
\end{align*}
for some constant $C'$. Passing to the limit as $r \to 0^+$ leads to
a contradiction, owing to assumption \eqref{ott21}. Equation \eqref{feb11} is thus established.
\\
Since the functions
$$
r \mapsto {\rm ess\, inf}_{B_r(x_0)} u\quad \hbox{and} \quad 
r \mapsto {\rm ess\, sup}_{B_r(x_0)}
u$$ are monotone in $r$,  they admit (finite) limits as $r \to 0^+$,
which, owing to
\eqref{feb11}, 
agree in $\Omega \setminus E$
In particular, the representative  $\widehat u$, defined by \eqref{utilde}, satisfies the equality
\begin{equation}\label{feb33}
\widehat u (x) = \lim _{r \to 0^+} {\rm ess\, inf}_{B_r(x)} u
= \lim _{r \to 0^+} {\rm ess\, sup}_{B_r(x)} u \quad \hbox{for every $x \in \Omega \setminus E$.}
\end{equation}
It is easily verified that the function 
$$
x \mapsto \lim _{r \to 0^+}  {\rm ess\, inf}_{B_r(x)} u$$
is lower-semicontinuous in $\Omega$. Hence, by \eqref{feb33}, $\widehat u$ is lower-semicontinuous  in $\Omega \setminus E$. Similarly, the function  $$
x \mapsto \lim _{r \to 0^+} {\rm ess\, sup}_{B_r(x)} u$$
is upper-semicontinuous in $\Omega$, and, by \eqref{feb33} again, $\widehat u$ is also upper-semicontinuos  in $\Omega \setminus E$. Altogether, we have shown that $\widehat u$ is continuous  in $\Omega \setminus E$. The proof is complete.
\qed
 
 \smallskip
 \par
 As mentioned in Section \ref{main}, Corollary \ref{orliczcapacitycor} follows from Theorem \ref{orliczcapacity} via the next result.
 
  \smallskip

\begin{proposition}\label{delta2lem}
Let $n \geq 2$, let $A$ be a Young function, and let $A_{n-1}$ be the Young function defined by \eqref{Apoincarebis}. If  $A \in \Delta _2$ near
infinity, then $A_{n-1}\in \Delta _2$ near infinity as well.
\end{proposition}
\par\noindent
{\bf Proof}. Owing to \eqref{deltanabla}, one has that
$\widetilde A \in \nabla _2$ near infinity. Thus, by \eqref{nabla2incr},   there exists
$\varepsilon_0>0$ such that, if $0<\varepsilon <\varepsilon_0$, 
the function $\widetilde A (t)t^{-1-\varepsilon}$ increasing for large $t$.
Hence, if $0< \varepsilon < \min\{\varepsilon_0, \tfrac 1{n-2}\}$,
then
\begin{align}\label{feb17}
\int _t^\infty \frac {\widetilde A(r)}{r^{1 + \frac {n-1}{n-2}}}\,
dr & \geq \frac {\widetilde A(t)}{t^{1 + \varepsilon}} \int _t
^\infty \frac {dr}{r^{\frac {n-1}{n-2}-\varepsilon}}  = \frac
{\widetilde A(t)}{(\frac {n-1}{n-2}-1-\varepsilon ) t^{ \frac
{n-1}{n-2}}} \quad \hbox{for large $t$.}
\end{align}
As a consequence,
\begin{align}\label{feb19}
\liminf _{t \to \infty} \frac{t
\frac{d \widetilde{A_{n-1}}}{dt}(t)}{\widetilde{A_{n-1}}(t)}& = \liminf _{t \to
\infty} \frac{t \Big( \frac{n-1}{n-2} t^{\frac {n-1}{n-2}-1}\int
_t^\infty \frac {\widetilde A(r)}{r^{1 + \frac {n-1}{n-2}}}\, dr -
\frac{\widetilde A (t)}t\Big)}{ t^{\frac {n-1}{n-2}} \int _t^\infty
\frac {\widetilde A(r)}{r^{1 + \frac {n-1}{n-2}}}\, dr}
\\ \nonumber & = \frac{n-1}{n-2} - \limsup _{t \to \infty}
\frac{\widetilde A (t)}{t^{\frac {n-1}{n-2}} \int _t^\infty \frac
{\widetilde A(r)}{r^{1 + \frac {n-1}{n-2}}}\, dr} \geq 1 +
\varepsilon \,.
\end{align}
By \eqref{feb19}, for every $\ep _1 \in (0, \ep)$ the function $\widetilde{A_{n-1}}(t) t^{-1-\ep_1}$ is increasing for large $t$. Hence, by \eqref{nabla2incr},
$\widetilde{A_{n-1}} \in \nabla _2$ near infinity, whence, thanks to \eqref{deltanabla},
namely $A_{n-1} \in \Delta _2$ near infinity. \qed

 \smallskip
\par\noindent

The link between the generalized capacities $C_{\Psi,1}$ and the classical Hausdorff measures $\mathcal H^{h(\cdot)}$ is discussed in the old paper \cite{taylor}. It provides us with a key tool in deriving Corollary  \ref{continuityhausdorff} from Theorem \ref{orliczcapacity}.

\smallskip
\par\noindent
{\bf Proof of Corollary \ref{continuityhausdorff}.}  Let $E$ be the set where the representative $\widehat u$ of $u$, exhibited in the proof of Theorem \ref{orliczcapacity}, is not continuous. Being the complement in $\Omega$ of the set where the limit of the averages of $u$ exists, the set $E$ is Borel measurable. Thus,  for every $k \in \mathbb N$, the set $E_k = E \cap B_k(0)$ is a bounded Borel set.
%
%
Given $k \in \mathbb N$,  let $f \in L^1_+(\rn
)$ be such that $I_\Psi f (x) \geq 1$ for $x \in E_k$, and let $\mu$
be any Borel measure, supported in $E_k$, such that $\mu (E_k)=1$. Then
\begin{align}\label{feb20}
1= \mu (E_k) \leq \int _{\rn}I_\Psi f (x)\, d \mu (x) = \int _{\rn}
I_\Psi \mu (y) f (y) \,dy \leq \|I_\Psi \mu \|_{L^\infty (\rn)}
\|f\|_{L^1(\rn )}.
\end{align}
Hence, from the very definition of $C_{\Psi,1}$-capacity,
\begin{align}\label{feb21}
C_{\Psi,1} (E_k) \geq \frac{1}{\|I_\Psi \mu \|_{L^\infty (\rn)}}.
\end{align}
Since, by Theorem \ref{orliczcapacity}, $C_{\Psi,1} (E_k) =0$ for every $k \in \mathbb N$, one has that $\|I_\Psi \mu \|_{L^\infty (\rn)}= \infty$  for every $\mu$ as above. This piece of information, combined with assumption \eqref{hB},   implies, via   \cite[Theorem
2]{taylor}, that $$\mathcal H ^{h(\cdot)} (E_k)
=0$$  for every  {$k\in \mathcal \N$}. Next, recall that, although   $\mathcal H^{h(\cdot)}$ is just  an outer measure, it is a measure when restricted to the class of Borel sets. Thus, since $E_k$ is an increasing sequence of Borel sets such that $E= \cup _{k=1}^\infty E_k$,
$$\mathcal H^{h(\cdot)} (E) = \lim _{k\to \infty}\mathcal H^{h(\cdot)} (E_k) = 0\,.$$ \qed
%
%
%
 
\smallskip
\par
Our last proof concerns Theorem \ref{febbraio0}.

\smallskip
\par\noindent
{\bf Proof of Theorem \ref{febbraio0}.}
%
Assume, the time being, that $u$ fulfills condition \eqref{finiteint}. It is clear from the proof of inequality \eqref{orlicz4} that it also holds with $u$ replaced by every approximating function $u_k$, defined as in the proof of Theorem \ref{orliczbound}.
Integrating over $(r, R)$ the inequality obtained from this replacement
 yields
\begin{equation}\label{orlicz5}
\int _r ^{R} t^{n-1} {A_{n-1}}\bigg(\frac {|u_k(x) -
u_k(y)|}{Ct}\bigg)\, dt \leq \int _{B_r(x_0)} A(|\nabla u_k|) \, dz
< \infty\,
\end{equation}
provided that $B_R(x_0) \subset \subset   \Omega$, $0<r<t<R$ and $x,y \in B_r(x_0)$.
  Assumption \eqref{ours} is equivalent to 
\begin{equation}\label{orlicz6}
\int _0 t^{n-1} A_{n-1}\Big(\frac 1t\Big)\, dt = \infty.
\end{equation}
As a consequence  of  \eqref{orlicz5} and \eqref{orlicz6}, the  sequence $u_k$ is  equi-continuous in $\Omega$.
 Thus, by Ascoli-Arzel\'a's theorem, the sequence $\{u_k\}$ converges
 to a continuous function $\overline u$, which agrees with $u$ a.e. in $\Omega$.
 Denote by $\omega (r)$ the modulus of continuity of $\overline u$ in $B_r(x_0)$.
 Inequality \eqref{orlicz5}
implies, via Fatou's lemma and equation \eqref{approx2}, that
\begin{equation}\label{modulus2}
\int _r ^{R}t^{n-1} A_{n-1}\bigg(\frac {\omega  (r)}{Ct}\bigg)\,
dt \leq \int _{B_R(x_0)} A(|\nabla u|)\, dz \,.
\end{equation}
From an application of inequality \eqref{modulus2} with $r= \tfrac R2$, and property \eqref{kt} applied with $A$ replaced by $A_{n-1}$, we deduce that, if $\omega (\tfrac R2) \geq 1$, then 
\begin{equation}\label{gen3}
\omega  (\tfrac R2) \int _{R/2} ^{R}t^{n-1} A_{n-1}\bigg(\frac {1}{Ct}\bigg)\,
dt  \leq \int _{R/2} ^{R}t^{n-1} A_{n-1}\bigg(\frac {\omega  (\tfrac R2)}{Ct}\bigg)\,
dt \leq \int _{B_R(x_0)} A(|\nabla u|)\, dz\,.
\end{equation}
Hence, 
\begin{equation}\label{gen4}
\omega  (\tfrac R2) \leq \max \Bigg\{1, \frac{\int _{B_R(x_0)} A(|\nabla u|) \,dz}{\int _{R/2} ^{R}t^{n-1} A_{n-1}(\frac {1}{Ct})\,
dt }\Bigg\}\,.
\end{equation}
Thus, if $0<r<\tfrac R2$, 
\begin{align}\label{gen5}
\int  _R^{\infty} t^{n-1}A_{n-1}\bigg(\frac {\omega 
(r)}{Ct}\bigg)\, dt
& = \bigg(\frac{\omega  (r)}{C}\bigg)^n\int _0 ^{\frac {\omega 
(r)}{RC}} \frac{A_{n-1}(s)}{s^{n+1}}\, ds \\ \nonumber & \leq \bigg(\frac{\omega  (\tfrac R2)}{C}\bigg)^n\int _0 ^{\frac {\omega
(\frac R2)}{RC}} \frac{A_{n-1}(s)}{s^{n+1}}\, ds < \infty\,,
\end{align}
where the equality holds by a change of variable in the integral, and the last inequality  by equation \eqref{febbraio4}.
Combining inequalities \eqref{modulus2} and \eqref{gen5} tells us that
\begin{align}\label{gen6}
\int _r ^{\infty} t^{n-1} {A_{n-1}}\bigg(\frac {\omega 
(r)}{Ct}\bigg)\, dt & = \int _r ^{R} t^{n-1} {A_{n-1}}\bigg(\frac {\omega 
(r)}{Ct}\bigg)\, dt + \int _R ^{\infty} t^{n-1} {A_{n-1}}\bigg(\frac {\omega 
(r)}{Ct}\bigg)\, dt \\ \nonumber & \leq  \int _{B_R(x_0)} A(|\nabla u|) \,dz + \bigg(\frac{\omega  (\tfrac R2)}{C}\bigg)^n\int _0 ^{\frac {\omega
(\frac R2)}{RC}} \frac{A_{n-1}(s)}{s^{n+1}}\, ds 
\end{align}
if $0<r\leq\tfrac R2$.
Denote by $C'=C'(x_0, R, u)$ the quantity on the rightmost side of  \eqref{gen6}, and set $C'' = \max \{C, C'\}$, where $C$ is the constant appearing in \eqref{gen6}. Thus,
\begin{equation}\label{modulus5}
\bigg(\frac{\omega  (r)}{C''}\bigg)^n\int _0 ^{\frac {\omega 
(r)}{rC''}} \frac{A_{n-1}(s)}{s^{n+1}}\, ds   = \int  _r^{\infty} t^{n-1}A_{n-1}\bigg(\frac {\omega 
(r)}{C''t}\bigg)\, dt \leq C'' \quad 
\hbox{if $0<r\leq\tfrac R2$.}
\end{equation}
or, equivalently,
\begin{equation}\label{modulus6}
B\bigg(\frac {\omega  (r)}{rC''} \bigg) \leq \frac {C''}{r^n} \quad \hbox{if $0<r\leq\tfrac R2$,}
\end{equation}
where $B$ is given by \eqref{G}. Hence,
$$\omega (r) \leq C'' r B^{-1}(C'' r^{-n}) \leq C'' {\rm max}\{C'', 1\}\,
r\, B^{-1}(r^{-n}) \quad \hbox{if $0<r\leq\tfrac R2$,}$$ where the last inequality holds by property \eqref{kt-1} applied to the Young function  $B$. As a consequence, $u \in C^{\omega (\cdot)}_{\rm loc}(\Omega)$.
\\
If assumption \eqref{finiteint} is dropped, the same argument yields this conclusion, when applied with
 $u$ replaced by $\lambda u$ for a suitable constant $\lambda >0$.
 \qed

\color{black}

\section{
Examples
}\label{examp}

\subsection{A customary example}\label{cust}
Results in the spirit of Theorem \ref{orliczbound} on $L^\infty$-estimates for  weakly monotone Orlicz-Sobolev functions, as well as information on the capacity of the exceptional set and on its Hausdorff measure, such as that provided in Theorems \ref{orliczcapacity} and \ref{hausdorff}, respectively, seem to be missing in literature. 
 In Example \ref{powerlog} below, we illustrate these results in 
a model instance of the Orlicz-Sobolev spaces $W^1_{\rm loc}L^p\log ^\alpha (\Omega)$  built upon  Young functions of power-logarithmic type. Let us emphasize that, as already pointed out  in Remark \ref{improveManfredi}, even the classical result 	\eqref{p>n-1} on the capacity of the singular set of weakly monotone functions is improved by Theorem \ref{orliczcapacity}, when the latter is specialized to standard Sobolev spaces $W^{1,p}_{\rm loc}(\Omega)$ corresponding to plain power type Young functions. 
 \par For those exponents $p$ and $\alpha$ that guarantee the continuity  everywhere  of weakly monotone functions from $W^1_{\rm loc}L^p\log ^\alpha (\Omega)$, we recover in Example \ref{powerlog} the results of \cite[Chapter 6]{IM}. We also compare the modulus of continuity of a weakly monotone function in $W^1_{\rm loc}L^p\log ^\alpha (\Omega)$ given by Theorem \ref{main} with that of an arbitrary function from the same Orlicz-Sobolev space, and  derive some interesting conclusions.

\begin{example}\label{powerlog} {\rm Assume that $A(t)=t^p \log^\alpha (c+t)$, where $p>1$, $\alpha \in \mathbb R$, for some positive constant $c$ so large that $A$ is a Young function.  Let us denote by $W^1_{\rm loc}L^p\log ^\alpha (\Omega)$ the local Orlicz-Sobolev space associate with $A$. 
By Theorem \ref{orliczbound}, any weakly monotone function in $u \in W^1_{\rm loc}L^p\log ^\alpha (\Omega)$ is locally bounded, and differentiable a.e. in $\Omega$, in any of the following cases:
\begin{equation}\label{ex1}
\begin{cases}
\hbox{$p=n$ and $\alpha < -1$, }
\\ \hbox{$n-1 <p <n$ and $\alpha \in \mathbb R$,}
\\ \hbox{either $n \geq 3$, $p=n-1$ and $\alpha > n-2$, or $n =2$, $p=1$ and $\alpha \geq 0$\,.}
\end{cases}
\end{equation}
Moreover, denote by $E$ the exceptional set of discontinuity points of $u$. Then an application of Corollary \ref{orliczcapacity} and Theorem \ref{hausdorff} tells us what follows:
\\
(i) If $p=n$ and $\alpha < -1$, then 
$$C_{t^n \log^{\alpha +1}t, 1}(E)=0\, \quad \hbox{and} \quad \mathcal H^{\log^{-\gamma} (1/s)}(E)=0\quad \hbox{for every $\gamma >-\alpha$}.$$
\\
(ii)
If either  $n-1<p<n$ and $\alpha \in \mathbb R$, or $n=2$, $p=1$ and $\alpha \geq 0$, then 
$$C_{t^p \log^{\alpha +1}t, 1}(E)=0\, \quad \hbox{and} \quad \mathcal H^{s^{n-p}\log^{-\gamma} (1/s)}(E)=0 \quad \hbox{for every $\gamma >1-\alpha$}.$$
\\
(iii)
If $n \geq 3$, $p=n-1$ and $\alpha >n-2$, then 
$$C_{t^{n-1} \log^{\alpha +3-n}t, 1}(E)=0\, \quad \hbox{and} \quad  \mathcal H^{s^{n-p}\log^{-\gamma} (1/s)}(E)=0 \quad \hbox{for every $\gamma >n-2-\alpha$}.$$
On the other hand, from Corollaries \ref{maindelta2} and \ref{recover} one deduces that $u$ is everywhere continuous in $\Omega$ in any of the following cases:
\begin{equation}\label{ex3}
\begin{cases}
\hbox{$p>n$ and $\alpha \in \mathbb R$, }
\\ \hbox{$p=n$ and $\alpha \geq -1$.}
\end{cases}
\end{equation}
Furthermore, from Theorem \ref{main} one infers that:
\\ (i) If $p>n$ and $\alpha \in \mathbb R$, then 
\begin{equation}\label{ex6}
u \in C^{r^{1-\frac np}\log^{-\frac \alpha p} (1/r)}_{\rm loc}(\Omega).
\end{equation}
\\ (ii) If $p=n$ and $\alpha >-1$, then 
\begin{equation}\label{ex7}
u \in C^{\log^{-\frac {\alpha +1} n} (1/r)}_{\rm loc}(\Omega).
\end{equation}
\\ (iii) If $p=n$ and $\alpha =-1$, then 
\begin{equation}\label{ex8}
u \in C^{(\log \log) ^{-\frac {1} n} (1/r)}_{\rm loc}(\Omega).
\end{equation}
The content of equations \eqref{ex3}--\eqref{ex8} recovers results from \cite[Chapter 6]{IM}.
\\ Let us notice that, if $p>n$ and $\alpha \in \mathbb R$, then the modulus  of continuity of a weakly monotone function $u$ given by \eqref{ex6} coincides with that of any Orlicz-Sobolev function in $u \in W^1_{\rm loc}L^p\log ^\alpha (\Omega)$ -- see \cite{CianchiRand}[Equation (6.13)] -- and, in particular, with that provided by the classical embedding theorem by Morrey, if $\alpha =0$. Thus, being a weakly monotone  does not provide a function with a better modulus of continuity in these cases. 
\\ By contrast, if $p=n$ and $\alpha >n-1$, then any  Orlicz-Sobolev function in $u \in W^1_{\rm loc}L^n\log ^\alpha (\Omega)$ is still continuous, but one  has just that
\begin{equation}\label{ex9}
u \in C^{\log ^{-\frac{\alpha - n+1}n}(1/r)}_{\rm loc}(\Omega)\,.
\end{equation}
This is of course a weaker property than \eqref{ex7}.  Weak monotonicity thus turns out to improve the quality of the modulus of continuity of Orlicz-Sobolev functions in borderline situations, a phenomenon that cannot be appreciated in the less fine scale of standard Sobolev spaces.


}
\end{example}

\subsection{Augmenting the existing literature: a non-standard example.}\label{nonstandard}


 Since condition \eqref{incr} implies \eqref{indexcond},  Corollary \ref{recover} recovers the results of \cite{IM} and \cite{KKMOZ}. The objective of Proposition \ref{improve} below is to demonstrate that  assumption \eqref{ours} of Theorem \ref{febbraio0} is actually weaker than the pair of assumptions \eqref{Iwcond}--\eqref{incr}.
 This shows that  Theorem \ref{febbraio0} can be applied to deduce the continuity of weakly monotone Orlicz-Sobolev  functions in certain situations where the criterion of \cite{IM} and \cite{KKMOZ} fails.

\begin{proposition}\label{improve}
Let $n \geq 2$. Then there exist  Young functions $A$ for which condition \eqref{ours} holds, whereas  \eqref{incr} fails.
\end{proposition}

Let us point out that the Young functions that will be exhibited in the proof of Proposition \ref{improve} satisfy, in addition, the $\Delta_2$-condition. Thus, they fulfill condition \eqref{ours} in the equivalent form \eqref{condelta2} appearing in Corollary \ref{maindelta2}.

\medskip
\par\noindent
{\bf Proof of Proposition \ref{improve}}. Let $A$ be a piecewise affine Young function  of the form \eqref{A}. Thus, there exists an increasing sequence $\{t_k\}$ of nonnegative numbers $t_k$, with $t_0=0$ and 
\begin{equation}\label{tkinf}
\lim_{k\to \infty} t_k\,=\,\infty\,,
\end{equation} 
and an increasing sequence 
$\{m_k\}$ of nonnegative numbers $m_k$, such that
\begin{equation}\label{defmk}
a(t)\,=\,m_k\qquad \qquad {\rm for}\; t\in (t_k,\, t_{k+1})\,,
\end{equation}
for $k \in \mathbb N$.
Hence, 
\begin{equation}\label{Function}
A(t)\,=\,\sum_{h=0}^{k-1}\,m_{h}(t_{h+1}-t_{h})\,+\, m_k(t-t_k),\qquad \hbox{if $t \in  [t_k, t_{k+1})$,}
\end{equation}
for $k \in \mathbb N$.
The conclusion will follow if we  show that, given any $q>1$, the sequences $\{t_k\}$ and $\{m_k\}$ can be chosen in such a way that the function $A$ satisfies condition \eqref{ours}, and the function
\begin{equation}\label{notincr}
t \mapsto \frac{A(t)}{t^{q}}\qquad \hbox{is not increasing.}
\end{equation}
Condition \eqref{notincr} is fulfilled if there exists $\alpha \in (1, q)$ such that 
\begin{equation*}
\lim_{t\to t_k^-}a(t)\,t\,-\,\alpha\, A(t)\,=\,0\,
\end{equation*}
for $k \in\mathbb N$, 
namely
\begin{equation}\label{nocondition}
a(t_{k}^-)\,t_k\,-\,\alpha\, A(t_k)\,=\,0\,
\end{equation}
for $k \in\mathbb N$.
 The following formulas can be derived 
from equations 
  \eqref{Function} and \eqref{nocondition}, via an induction argument:
\begin{align}\label{2}
t_{k+1}\,&=\, \frac{\alpha t_1(m_1-m_0)(\alpha m_2-m_1)\cdots (\alpha m_{k-1}-m_{k-2}) (\alpha m_k-m_{k-1})}{(\alpha -1)^k\, m_1m_2\cdots m_k} \\ \nonumber  &=\,\frac{\alpha t_1}{(\alpha -1)^k}\,\Big(1-\frac{m_0}{m_1}\Big)\Big(\alpha-\frac{m_1}{m_2}\Big)\cdots \Big(\alpha-\frac{m_{k-2}}{m_{k-1}}\Big) \Big(\alpha-\frac{m_{k-1}}{m_k}\Big)\,,\end{align}
and
\begin{align}\label{2'}
t_{k+1}\,-\, t_k\, & =\, \frac{\alpha t_1(m_1-m_0)(\alpha m_2-m_1)\cdots (\alpha m_{k-1}-m_{k-2})(m_{k}-m_{k-1})}{(\alpha -1)^k\, m_1m_2\cdots m_k}
\\ \nonumber  &=\,\frac{\alpha t_1}{(\alpha -1)^k}\,\Big(1-\frac{m_0}{m_1}\Big)\Big(\alpha-\frac{m_1}{m_2}\Big)\cdots \Big(\alpha-\frac{m_{k-2}}{m_{k-1}}\Big) \Big(1-\frac{m_{k-1}}{m_{k}}\Big)\,
\end{align}
for $k \in\mathbb N$.
Let us define 
 the sequence $\{m_k\}$ in such a way that
\begin{equation}\label{mkdef}
\beta\, m_k\,=\, \alpha\,m_k\,-\,m_{k-1}\, 
\end{equation}
for $k \in \mathbb N$, with $\beta \in (\alpha -1, \alpha)$ to be fixed later. Hence,
\begin{equation}\label{mk}
m_k\,= \,\frac{m_{0}}{(\alpha-\beta)^k}
\end{equation}
and 
\begin{equation}\label{mkbis}
m_k\,-\,m_{k-1}\,=\, (\beta-\alpha +1) \frac{m_0}{(\alpha-\beta)^k}\,
\end{equation}
for $k \in\mathbb N$.
By \eqref{2},
\begin{align}\label{1}
t_{k+1} 
= \frac{\alpha t_1 (m_1-m_0)\beta^{k-1}}{m_1(\alpha -1)^k}\,
\end{align}
for $k \in \mathbb N$, whence \eqref{tkinf} holds.  Furthermore, by \eqref{2'},
\begin{align}\label{1'}
t_{k+1} -t_k =  \frac{\alpha t_1 (m_1-m_0)(\beta -\alpha +1)\beta^{k-2}}{m_1(\alpha -1)^k}\,
\end{align}
for $k \in \mathbb N$.
\\
Assume first that $n=2$. We claim that choosing $\alpha \in (1, \min\{q, 2\})$ and $\beta =1$ in \eqref{mkdef} yields a function $A$ fulfilling condition \eqref{ours}, which reads
\begin{equation}\label{nov1}
\int^{\infty}_{t_1}\frac{A(t)}{t^3}\,dt = \infty\,
\end{equation}
for $n=2$.
Indeed, owing to property \eqref{gen9}, equation \eqref{nov1} is equivalent to 
\begin{equation}\label{nov1bis}
\int_{t_1}^{\infty}\frac{a(t)}{t^2}\,dt = \infty\,.
\end{equation}
On the other hand, by \eqref{mk}, \eqref{1} and \eqref{1'}, with $\beta =1$,
\begin{align}\label{3}
\int^{\infty}_{t_1}\frac{a(t)}{t^2}\,dt =
 \sum^{\infty}_{k=1}\int_{t_k}^{t_{k+1}}\frac{m_k}{t^2}\,dt =
\frac{m_0m_1}{t_1(m_1-m_0)}\frac{2-\alpha}{\alpha-1}\sum^{\infty}_{k=1}1 = \infty\,.
\end{align}
whence \eqref{nov1bis} follows.
\\ Assume next that $n \geq 3$. Let us preliminarily observe that any function $A$ defined by \eqref{Function}, with $m_k$ obeying \eqref{mk}, satisfies the $\Delta_2$-condition. Owing to the equivalence of equations \eqref{delta2} and \eqref{delta2'}, in order to verity this assertion it suffices to show that  there exists a constant $c$ such that
\begin{equation}\label{exdelta1}
a(2t)\le c\, a(t) \qquad \hbox{for $t \geq0$.}
\end{equation}
Given  $t>0$, let $k\in \N$  be the index satisfying 
$t\,\in\,[t_k,\, t_{k+1})$,
whence $2t\,\in\,[2t_k,\, 2t_{k+1})$. If we prove that there exists a constant $c$ such that
\begin{equation}\label{exdelta5}
a(2t_{k+1})\le c\, a(t_k)\,
\end{equation}
for  $k\in\N$, 
then inequality \eqref{exdelta1} will follow, inasmuch as
\begin{equation*}
a(2t)\le a(2t_{k+1})\le c\, a(t_k)\le c\, a(t)\, \qquad \hbox{for $t>0$.}
\end{equation*}
Given $k \in \mathbb N$, denote by  $j=j(k)\in \N$  the index fulfilling
\begin{equation}\label{exdelta2}
2t_{k+1}\,\in\,[t_{j},\, t_{j+1})\,.
\end{equation}
Thanks to equation \eqref{1}, condition \eqref{exdelta2} is equivalent  to the inequalities
\begin{equation}\label{nov3}
\frac{\alpha t_1 (m_1-m_0)}{m_1(\alpha-1)}\left(\frac{\beta}{\alpha-1}\right)^{j-2}\le 
\frac{2\alpha t_1 (m_1-m_0)}{m_1(\alpha-1)}\left(\frac{\beta}{\alpha-1}\right)^{k-1}< \frac{\alpha t_1 (m_1-m_0)}{m_1(\alpha-1)}\left(\frac{\beta}{\alpha-1}\right)^{j-1}\,.
\end{equation}
On setting $b\,=\,\frac{\beta}{\alpha-1}$, equation \eqref{nov3} is in turn equivalent to
$$ j-2\,\le\,\log_b2+k-1\,\le\,j-1\,,$$
whence
$$ j\,\le\,1+\log_b2+k\,\le\,j+1\,\le\, [\gamma]+k+2$$
where $\gamma =1+\log_b2$. Therefore
\begin{align}
a(2t_{k+1})\le \,a(t_{j+1})& \le \, a(t_{[\gamma]+k+2})\,=\,m_{[\gamma]+k+2}  \\ \nonumber  &=\,
\frac{m_0}{(\alpha -\beta)^{[\gamma]+k+2}}  =\,\frac{m_k}{(\alpha -\beta)^{[\gamma]+2}}\,=\, \frac{1}{(\alpha -\beta)^{[\gamma]+2}}a(t_k)
\end{align}
We conclude that inequality \eqref{exdelta5}, and hence \eqref{exdelta1}, holds with $c\,=\,\frac{1}{(\alpha -\beta)^{[\gamma]+2}}$.
\\ Since $A$ satisfies the $\Delta _2$-condition, by \cite[Lemma 3.3]{CarozzaCianchi} condition \eqref{ours} is equivalent to \eqref{condelta2}. By property \eqref{gen9},  the latter is in turn equivalent to 
\begin{equation}\label{ex2}
 \int ^{\infty}\frac{1}{a(t)^\frac{2}{n-2}}\bigg(\int_t^{\infty}\frac{1}{a(s)^{\frac{1}{n-2}}}\,ds\bigg)^{-n}\,dt =\, \infty\,.
\end{equation}
We shall show that \eqref{ex2} holds provided that $\beta$ is chosen in such a way that $\alpha - \beta$ is sufficiently small. Indeed, if $k \in \mathbb N$ and $t \in [t_k, t_{k+1}]$, then, by \eqref{mk} and \eqref{1'}, 
\begin{align}\label{nov4}
\int_t^{\infty}\frac{1}{a(s)^{\frac{1}{n-2}}}\,ds\,& \le\, \int_{t_k}^{\infty}\frac{1}{a(s)^{
\frac{1}{n-2}}}\,ds\,= \sum _{h=k}^\infty \int_{t_h}^{t_{h+1}}\frac{1}{a(s)^{
\frac{1}{n-2}}}\,ds
\\ \nonumber & =
\,\sum_{h=k}^{\infty}\frac{t_{h+1}-t_h}{m_h^{\frac{1}{n-2}}}\,=\,C\sum_{h=k}^{\infty}\bigg[\frac{\beta(\alpha-\beta)^{\frac{1}{n-2}}}{\alpha-1}\bigg]^{h}
\,=\,C'\bigg[\frac{\beta(\alpha-\beta)^{\frac{1}{n-2}}}{\alpha-1}\bigg]^{k}
\end{align}
for suitable constants $C$ and $C'$, provided that $\alpha - \beta$ is sufficiently small.
From inequality \eqref{nov4} and equations \eqref{mk} and \eqref{1'} one then infers that
\begin{align}\label{nov5}
\int_{t_1} ^{\infty} &\frac{1}{a(t)^\frac{2}{n-2}}\bigg(\int_t^{\infty}\frac{1}{a(s)^{\frac{1}{n-2}}}\,ds\bigg)^{-n}\,dt \,= \sum_{k=1}^{\infty}
\int_{t_k} ^{t_{k+1}}\frac{1}{a(t)^\frac{2}{n-2}}\,dt\bigg(\int_t^{\infty}\frac{1}{a(s)^{\frac{1}{n-2}}}\,ds\bigg)^{-n} dt \\ \nonumber
& \ge  C\sum_{k=1}^{\infty}\frac{t_{k+1}-t_k}{m_k^{\frac{2}{n-2}}} 
\bigg[\frac{\beta(\alpha-\beta)^{\frac{1}{n-2}}}{\alpha-1}\bigg]^{-nk}
 = C'\sum_{k=1}^\infty \bigg[\bigg(\frac{ \alpha -1}{\beta}\bigg)^{n-1}\frac 1{\alpha-\beta}\bigg]^{k}
\end{align}
for suitable positive constants $C$ and $C'$. Since the last series diverges if $\alpha - \beta$ is small enough, equation \eqref{ex2} follows. \qed

\medskip
 \noindent {\bf Acknowledgments.} This research was
partly supported by the Research Project of Italian Ministry of
University and Research (MIUR)  2012TC7588 ``Elliptic and parabolic
partial differential equations: geometric aspects, related
inequalities, and applications" 2012,  and by  GNAMPA of Italian INdAM
(National Institute of High Mathematics).


\begin{thebibliography}{99}

\bibitem[AH]{AH}  D.R.Adams \& L.I.Hedberg, \lq\lq Function spaces and potential theory", Springer, Berlin (1996).

\bibitem[AC]{AlCi} A.Alberico \& A.Cianchi,
 Differentiability properties of Orlicz-Sobolev functions,
  \emph{Ark. Math.} \textbf{43} (2005),  1--28.


\bibitem[ALT]{ALT} A.Alvino, P.-L.Lions \& G.Trombetti, On optimization problems with prescribed rearrangements,
  \emph{Nonlinear Anal.} \textbf{13 } (1989),   185--220.



\bibitem[Ba]{B} J.Ball, Convexity conditions and existence theorems in nonlinear elasticity, \emph{ Arch. Rat. Mech. Anal.} \textbf{63} (1977), 337--403.



\color{black}
\bibitem[CC]{CarozzaCianchi} M.Carozza  \& A.Cianchi,  Smooth approximation of Orlicz-Sobolev maps between Riemannian manifolds,
 \emph{Potential Analysis} \textbf{45}
(2016), 557--578.

\bibitem[Ci1]{Cianchi_ASNS} A.Cianchi,
 Continuity properties of functions from Orlicz-Sobolev spaces and
  embedding theorems, \emph{Ann. Scuola Norm. Sup. Pisa}
 \textbf{23} (1996), 575--608.

\bibitem[Ci2]{Cianchiibero}  A.Cianchi,  Optimal Orlicz-Sobolev embeddings,
 \emph{Rev. Mat. Iberoamericana} \textbf{20}
(2004), 427--474.


\bibitem[CM]{CianchiMusil}  A.Cianchi \& V.Musil,  Optimal  domain spaces in Orlicz-Sobolev  embeddings,
 \emph{Indiana Univ. Math. J.}, to appear.

\bibitem[CR]{CianchiRand}  A.Cianchi \& M.Randolfi,  
On the modulus of continuity of weakly differentiable functions, \emph{Indiana Univ. Math. J.} {\bf 60} (2011),   1939--1973. 


\bibitem[FM]{FM}  T.Futamura \& Y.Mizuta,   Continuity of weakly monotone Sobolev functions,
 \emph{Adv. Math. Sci. Appl.} \textbf{15}
(2005), 571--585.



\bibitem[GI]{GI} F.W.Gehring \& T.Iwaniec, The limit of mappings with finite distortion,
 \emph{Ann. Acad. Sci. Fenn. Math. } \textbf{24}
(1999), 253--264.

\bibitem[GV]{GV} V.M.Gol'dstein \&  S.K.Vodop'yanov, Quasiconformal mappings and spaces of
functions with generalized first derivatives,
 \emph{Sibirsk. Mat. Z. } \textbf{17}
(1976), 515--531.

\bibitem[HK]{HK} J.Heinonen \&  P.Koskela, Sobolev mappings with integrable dilatations,  \emph{Arch. Rational Mech. Anal.} \textbf{125} (1993) 81--97.

\bibitem[HKM]{HKM} J.Heinonen, T.Kilpelainen \& O.Martio, \lq\lq Nonlinear potential theory of
degenerate elliptic equations", Oxford Univ. Press, 1993.

\bibitem[IM]{IM} T.Iwaniec \& G.Martin, \lq\lq Geometric function theory and nonlinear analysis", Clarendon
Press, Oxford, 2001.

\bibitem[IKO]{IKO} T.Iwaniec, P.Koskela \& J.Onninen, Mappings of finite distortion:
monotonicity and continuity,
 \emph{Invent. math.} \textbf{144}
(2001), 507--531.

\bibitem [IS]{IS} T.Iwaniec \& V.Sverak, On Mappings with integrable dilatation,
 \emph{Proc. Amer. Math. Soc.} \textbf{118}, (1993), 181--188.
 
 \bibitem [KKM]{KKM} J.Kauhanen, P.Koskela \& J.Maly,
   On functions with derivatives in a Lorentz space, \emph{Manuscripta Math.} \textbf{100} (1999), 87--101.
 
 
\bibitem [KKMOZ]{KKMOZ} J.Kauhanen, P.Koskela, J.Maly,
J.Onninen \& X.Zhong, Mappings of finite distortion:
sharp Orlicz-conditions,
 \emph{Rev. Mat. Iberoamer.} \textbf{19}
 (2003), 857--872.
 
\bibitem[KMV]{KMV} P.Koskela, J.J.Manfredi \& E.Villamor,  Regularity theory and traces of
A-harmonic functions,
 \emph{Trans. Amer. Math. Soc.} \textbf{348}
(1996), 755--766.

\bibitem[KO]{KO} P.Koskela \& J.Onninen,  Mappings of finite distorsion:
the sharp modulus of continuity,
 \emph{Trans. Amer. Math. Soc.} \textbf{355} (2003), 1905--1920.
 
 \bibitem[KR]{KR} M.A.Krasnosel'skii \& Ja.B.Rutickii, ``Convex functions and Orlicz spaces", P. Noordhoff Ltd., Groningen, 1961.
 
\bibitem[Le]{L} H.Lebesgue, Sur le probleme de Dirichlet,
 \emph{Rend. Cir. Mat. Palermo} \textbf{24}
(1907), 371--402.

\bibitem [Man]{M} J.J.Manfredi, Weakly monotone functions,
 \emph{J. Geom. Anal.} \textbf{4}
 (1994),  393--402.
 
 \bibitem [MV1]{MV1} J.J.Manfredi \& E.Villamor, Traces of monotone Sobolev functions,
 \emph{J. Geom. Anal.} \textbf{6}
 (1996),  433--434.
 
\bibitem [MV2]{MV2} J.J.Manfredi \& E.Villamor, Traces of monotone Sobolev functions in
weighted Sobolev spaces,
 \emph{Illinois J. Math.} \textbf{45}
 (2001),  403--422.
 
 \bibitem[Maz]{Mazya} V.G.Maz'ya, \lq\lq Sobolev spaces with applications to partial differential equations",  Springer-Verlag, Berlin, 2011.
 
  \bibitem[O'N]{oneil}
  R.O'Neil,   Convolution operators in $L(p,q)$
spaces, \emph{Duke Math. J.} {\bf 30} (1963), 129--142.


\bibitem [On]{O} J.Onninen, Differentiability of monotone Sobolev functions, \emph{Real Anal. Exchange} 
\textbf{26} (2000), 761--772.

\bibitem [OZ]{OZ} J.Onninen \& X.Zohng, A note on mappings of finite distortion: the sharp modulus of continuity, \emph{Michigan Math. J.} 
\textbf{53} (2005), 329--335.




 \bibitem[RR1]{RR1}  M.M.Rao \& Z.D.Ren, ``Theory of Orlicz spaces", Marcel Dekker Inc.,
 New York, 1991.
 
 \bibitem[RR2]{RR2}  M.M.Rao \& Z.D.Ren, ``Applications of Orlicz spaces", Marcel Dekker Inc.,
 New York, 2002.
 
\bibitem[Re]{R} Yu.G.Reshetnyak, ``Space mappings with bounded distortion", Trans. Math.
Monographs Amer. Math. Soc. 73, 1989.

 
  \bibitem [Ste]{Step} V.Stepanoff, Sur les conditions de la differentielle totale, \emph{Mat. Sb.} {\bf 32} (1925), 511-526.
 
 \bibitem [Str]{Strom} J.-O.Str\"omberg,
Bounded mean oscillation with {O}rlicz norms and duality of {H}ardy
  spaces,
 {\em Indiana Univ. Math. J.}  {\bf 28} (1979), 511--544.
 
 
 
 
\bibitem [Sv]{S} V.Sverak, Regularity properties of deformations with finite energy,
 \emph{Arch. Rat. Mech. Anal.} \textbf{100} (1988) 105--127.



\bibitem[Tal]{Talenti} G.Talenti, An embedding theorem, in \emph{Essays of Math. Analysis in honour
of E. De Giorgi,} Birkhauser Verlag, Boston, 1989.

\bibitem [Tay]{taylor}  S.J.Taylor,   On the connexion between Hausdorff measures and generalized capacity,
 \emph{Proc. Cambridge Philos. Soc.} \textbf{57} (1961) 524--531.
 

\end{thebibliography}
\end{document}